# Multi-Patch Isogeometric Convolution Hierarchical Deep-learning Neural Network

Lei Zhang[a,b,1], Chanwook Park[c,1], T.J.R. Hughes[d], Wing Kam Liu[c†]

[a]The State Key Laboratory of Nonlinear Mechanics, Institute of Mechanics, Chinese Academy of Sciences, Beijing 100190, China
[b]School of Engineering Science, University of Chinese Academy of Sciences, Beijing 100049, China
[c]Department of Mechanical Engineering, Northwestern University, Evanston 60208, USA
[d]Institute for Computational Engineering and Sciences, The University of Texas at Austin, 201 East 24th Street, 1 University Station C0200, Austin, TX 78712, USA
(1): These authors contributed equally to the manuscript.
(†): Corresponding author, email: w-liu@northwestern.edu

**Multi-Patch Convolution – IsoGeometric Analysis (C-IGA) Highlights**

- C-IGA formulation on multi-Computer-Aided Design (CAD)-patch systems is developed.
- Given the coarsest CAD mapping, one can refine the mesh in the physical domain arbitrarily using any meshing tools.
- Two compatibility conditions are mathematically and numerically proven; 1) *nodal compatibility* where the physical nodes on a patch interface are matched, and 2) $G^0$ *compatibility* where the entire interpolation field across the interface is matched.

**Abstract**

A seamless integration of neural networks with Isogeometric Analysis (IGA) was first introduced in [1] under the name of Hierarchical Deep-learning Neural Network (HiDeNN) and has systematically evolved into Isogeometric Convolution HiDeNN (in short, C-IGA) [2]. C-IGA achieves higher order approximations without increasing the degree of freedom. Due to the Kronecker delta property of C-IGA shape functions, one can refine the mesh in the physical domain like standard finite element method (FEM) while maintaining the exact geometrical mapping of IGA. In this article, C-IGA theory is generalized for multi-CAD-patch systems with a mathematical investigation of the compatibility conditions at patch interfaces and convergence of error estimates. Two compatibility conditions (nodal compatibility and $G^0$ (i.e., global $C^0$) compatibility) are presented and validated through numerical examples.

## 1. Introduction

The universal approximation theorem [3] that a deep neural network (DNN) can represent any functional relationships inspired the usage of neural networks as an approximate solution of a partial differential equation (PDE) [4]. Early and highly cited works [5-7] directly adopt neural networks as the approximate solution and exploit advanced programming tools such as TensorFlow [8] or PyTorch [9] with auto-differentiation to determine the solution.

Physics-informed neural network (PINN) [5] is the most widely used terminology to refer to these works. PINN consists of three steps: 1) approximate field variables with DNN; 2) at randomly spread collocation points, measure the residual of the strong form (or weak form if it is variational PINN [10]), boundary condition, and initial condition to formulate a loss function; and 3) optimize weights and biases of DNN by minimizing the loss function. The advantages of PINNs include a straightforward formulation of the solution field and loss function, as well as the ease of implementation using cutting-edge programming tools (e.g., [6] state that PINN has benefits over conventional numerical methods for high dimensional problems or problems with complex geometry where meshing is arduous or impossible). In contrast, the downside of PINNs is its inefficient solution scheme and inaccuracy compared to classical numerical methods (e.g., finite element method (FEM)), especially for low dimensional and regular-shaped problems [11, 12]. The training process of PINN where the strong form, boundary condition, and initial condition are imposed is the major bottleneck toward fast convergence and accurate solutions [11].

As opposed to PINN, several studies have been reported that link DNN with the FEM [1, 13-16]. These approaches utilize partially (or fully) connected DNN to represent the finite element interpolation functions. In [1], the authors propose a hierarchical deep-learning neural network – FEM (HiDeNN-FEM) that replaces the whole FEM process with partially connected hierarchical DNNs. As HiDeNN shape functions are mathematically equivalent to those used in FEM, the solution accuracy of HiDeNN is comparable to FEM. HiDeNN can reproduce any higher order shape functions in FEM, meshfree, and IGA, resulting in a faster convergence of error estimate [1]. Additionally, these neural networks are fully interpretable, that is, every weight and bias can be expressed as a function of nodal coordinates. For this reason, one can update nodal coordinates during the solution process like r-mesh adaptivity [1, 15].

Convolution HiDeNN (C-HiDeNN) takes one step further from HiDeNN [17, 18] to realize compact-supported (or non-local to the extreme) nonlinear interpolation. C-HiDeNN has one extra hidden layer to the HiDeNN shape function whose weights are called *convolution patch functions*. Instead of adding extra nodes for each element (like higher order Lagrange elements in FEM), C-HiDeNN leverages neighboring nodes to construct convolution patch functions whose nonlinearity is determined by the basis functions adopted in the convolution patch functions. For this reason, C-HiDeNN can use linear finite element meshes (e.g., 4-node tetrahedral element or 8-node brick element in 3D) but can achieve higher polynomial order and accuracy without increasing the degrees of freedom (DoFs) [17, 18]. In [17], the authors demonstrated faster convergence and numerical efficiency of C-HiDeNN running on graphics processing unit (GPU). C-HiDeNN theory was further extended to C-HiDeNN Tensor Decomposition (C-HiDeNN-TD) that can dramatically reduce the global DoFs for manufacturing simulations and topology optimization [18, 19].

Since C-HiDeNN(-TD) is formulated under a mesh-based geometrical representation, it has been limited in achieving a seamless integration of computer-aided design (CAD) models and numerical analysis (also known as computer-aided engineering (CAE)), which Isogeometric Analysis (IGA) is renowned for [20]. IGA replaces Lagrange polynomials in FEM shape function with spline functions used in CAD, enabling a direct usage of geometric basis functions for the solution approximation [20, 21].

The same philosophy – replacement of polynomial basis functions with splines – inspired the development of isogeometric C-HiDeNN (in short, C-IGA) [2]. Because the C-HiDeNN shape function can be constructed to reproduce any nonlinear function, one can replace the polynomial basis functions in the convolution patch functions with spline functions to reproduce the exact geometry of a CAD model. The mathematical proofs demonstrating the equivalence between the geometrical mappings employed in IGA and C-IGA are provided in [2]. C-IGA presents a new direction for the development of C-HiDeNN towards the compatible workflows between CAD and CAE.

So far, C-IGA has been developed only for a single B-spline or NURBS patch [2]. However, most CAD models for complex geometries are built with multiple patches [22, 23], emphasizing the need for multi-CAD-patch (or multi-patch for brevity) C-IGA. The biggest issue that arises with a multi-patch formulation is the compatibility: how the continuity of a solution field across patch interfaces is guaranteed. The goal of this paper is to develop multi-patch C-IGA theory that resolves this compatibility issue.

In this paper, the formulation, convergence theorem, and implementation of multi-patch Isogeometric **C**onvolution **H**ierarchical **De**ep-learning **N**eural **N**etwork (C-HiDeNN) (in short, multi-patch C-IGA) are investigated. It starts with a comprehensive review on multi-patch IGA theory and C-HiDeNN/C-IGA theory in section 2. Section 3 provides mathematical preliminaries of multi-patch C-IGA theory including interpolation estimates and convergence criteria. In Section 4, two fundamental theorems on multi-patch C-IGA are provided: one for the nodal compatibility and the other for the $G^0$ (i.e., globally $C^0$ continuous) compatibility. Numerical proofs on the convergence of multi-patch C-IGA then follows in Section 5.

Table 1 Nomenclature

| Symbol | Description |
|---|---|
| $\Omega, \hat{\Omega}, \tilde{\Omega}$ | Physical, parametric, parent domains, respectively. |
| $\Gamma, \hat{\Gamma}, \tilde{\Gamma}$ | Boundary of physical, parametric, parent domains, respectively. |
| $\boldsymbol{x}, \boldsymbol{\xi}, \tilde{\boldsymbol{\xi}}$ | Physical, parametric, parent coordinates, respectively. |
| $d$ | Dimension of input coordinates, $\boldsymbol{x}, \boldsymbol{\xi}, \tilde{\boldsymbol{\xi}} \in \mathbb{R}^d$ |
| $\hat{B}_{i,p}(\xi) / \hat{R}_{i,p}(\xi)$ | $i$-th single variate B-spline/NURBS basis of polynomial order $p$ |
| $p_m, n_m$ | Polynomial order and the number of basis functions of $m$-th dimension |
| $\mathcal{R}$ | Set of IGA indices. When $d = 2$, $\mathcal{R} = \{(i,j): i \in \{1,\dots, n_{m=1}\}, j \in \{1,\dots, n_{m=2}\}\}$ |
| $\hat{B}_{(i,j)}(\boldsymbol{\xi}) / \hat{R}_{(i,j)}(\boldsymbol{\xi})$ | $(i,j)$-th (double index) bivariate B-spline/NURBS basis when $d = 2$ |
| $w_{(i,j)}, \boldsymbol{P}_{(i,j)}, c_{(i,j)}$ | $(i,j)$-th (double index) control weights, points, and variables when $d = 2$ |
| $\boldsymbol{k}$ | Singe index. E.g., $(i,j) \to \boldsymbol{k}, \boldsymbol{k} \in \mathcal{R}$ |
| $\boldsymbol{F} : \hat{\Omega} \to \Omega$ | Forward IGA mapping from the parametric domain to the physical domain |
| $\boldsymbol{F}^{-1} : \Omega \to \hat{\Omega}$ | Inverse IGA mapping from the physical domain to the parametric domain |
| $B_{\boldsymbol{k}}(\boldsymbol{x}) / R_{\boldsymbol{k}}(\boldsymbol{x})$ | B-spline/NURBS basis with single index $\boldsymbol{k}$ in physical domain $\Omega$ |
| $V_h$ | Solution space of single-patch IGA |
| $\alpha, \beta$ | CAD patch index, mostly appears as a superscript of a variable |
| $\Gamma^{(\alpha,\beta)}$ | Interface of patch $\Omega^{(\alpha)}$ and $\Omega^{(\beta)}$ |
| $s, a, p$ | Patch size (non-negative integer), dilation parameter, polynomial order (non-negative integer). Hyper-parameters of C-HiDeNN/C-IGA |
| $A^e, A_s^{x_I}, A_s^e$ | Set of nodes in element $e$, in nodal convolution patch centered at node $x_I$ with convolution patch size $s$, and in elemental convolution patch of element $e$ with convolution patch size $s$. |
| $N_I(\tilde{\xi})$ | FEM shape function for node $x_I$ |
| $W_{s,a,p,K}^{x_I}/W_{s,a,p,K}^{\xi_I}$ | Convolution patch function of C-HiDeNN/C-IGA for node $x_K/\xi_K$ defined in the nodal convolution patch centered at node $x_I/\xi_I$ |
| $\tilde{N}_J(\tilde{\xi})/\tilde{N}_J(\xi)$ | C-HiDeNN/C-IGA shape function of node $J$ |
| $u_J$ | Nodal variable of node $J$ in C-HiDeNN/C-IGA |
| $u^h(x)/u^h(\xi)$ | C-HiDeNN/C-IGA approximate solution |
| $\mathcal{J}$ | General interpolation operator |
| $H^p(\Omega)/H_b^p(\Omega)$ | $p$-th order Sobolev/broken Sobolev space in the physical domain $\Omega$ |
| $\Gamma^D/\Gamma^N$ | Dirichlet/Neumann boundary of the physical domain $\Omega$ |
| $a(\cdot,\cdot)$ | Bilinear form of differential operator |
| $\langle u \rangle_{\Gamma^{(\alpha,\beta)}}/[u]_{\Gamma^{(\alpha,\beta)}}$ | Mean value / jump term at the patch interface $\Gamma^{(\alpha,\beta)}$ |
| $\mathcal{S}^h/\mathcal{V}^h$ | Trial / test function space |
| $\rho$ | Penalty coefficient |
| $u^{h,(\alpha)}(\xi)$ | C-IGA approximate solution at patch $\alpha$ |

## 2. Preliminaries
### 2.1. Isogeometric Analysis (IGA)
#### 2.1.1. Single-Patch IGA

To start with, a brief introduction to the B-splines and non-uniform rational B-splines (NURBS) basis functions will be provided. In a one-dimensional domain, we define positive integers: $p$ for the polynomial order of the basis function; and $n$ for the number of basis functions. Then an open knot vector is defined in a unit interval $\hat{\Omega} = [0,1]$ as $\Xi := [\xi_1, ..., \xi_{n+p+1}]$ where $\xi_i \leq \xi_{i+1}$, $\xi_1 = 0$, $\xi_{n+p+1} = 1$, the multiplicity of the first and the last knots are $p+1$, and that of the other knots is at most $p$. The B-spline basis function is defined based on the Cox-de Boor recursive formula [24]. Starting with $p = 0$,

$$\hat{B}_{i,0}(\xi) = \begin{cases} 1 & if\ \xi_i \leq \xi < \xi_{i+1}, \\ 0 & otherwise, \end{cases} \quad (1)$$

and for $p \geq 1$ the B-spline functions follow the recursion:

$$\hat{B}_{i,p}(\xi) = \frac{\xi - \xi_i}{\xi_{i+p} - \xi_i}\hat{B}_{i,p-1}(\xi) + \frac{\xi_{i+p+1} - \xi}{\xi_{i+p+1} - \xi_{i+1}}\hat{B}_{i+1,p-1}(\xi), \quad (2)$$

with an assumption of $0/0 = 0$. Since the open knot vector has $p+1$ repetitions at both ends, it can interpolate at the ends [20].

The multi-variate B-splines and NURBS can be defined through tensor product. For a $d$-dimensional problem (typically $d = 2,3$), the open knot vector of $m$-th dimension ($m = 1, ..., d$) with $p_m$ polynomial order and $n_m$ basis functions is expressed as $\Xi_m := [\xi_{m,1}, ..., \xi_{m,n_m+p_m+1}]$, where $\xi_{m,i}$ is the $i$-th knot of the $m$-th dimension. The B-splines for each dimension, $\hat{B}_{i,p_m}(\xi_m)$, are computed from Eq. ( 1 )-( 2 ). For the sake of simplicity, we set $d = 2$, a two-dimensional problem. The set of all double indices of B-splines or NURBS basis functions is denoted by:

$$\mathcal{R} = \{(i,j): i \in \{1, \ldots, n_{m=1}\}, j \in \{1, \ldots, n_{m=2}\}\}. \tag{3}$$

Through the tensor product, the bivariate B-splines are defined as

$$\hat{B}_{(i,j)}(\xi_1, \xi_2) = \hat{B}_{i,p_1}(\xi_1)\hat{B}_{j,p_2}(\xi_2), \tag{4}$$

and the bivariate NURBS basis functions are given by:

$$\hat{R}_{(i,j)}(\xi_1, \xi_2) = \frac{\hat{B}_{i,p_1}(\xi_1)\hat{B}_{j,p_2}(\xi_2)w_{(i,j)}}{W(\xi_1, \xi_2)}, \tag{5}$$

with the weighting function $W(\xi_1, \xi_2) = \sum_{(i,j) \in \mathcal{R}} \hat{B}_{i,p_1}(\xi_1)\hat{B}_{j,p_2}(\xi_2)w_{(i,j)}$,

where $w_{(i,j)}, (i,j) \in \mathcal{R}$ are the positive weights called _control weight_. For each double index $(i,j)$, a _control point_ $\boldsymbol{P}_{(i,j)} \in \mathbb{R}^2$ or $\boldsymbol{P}_{(i,j)} \in \mathbb{R}^3$ is defined in the physical domain $\Omega$ to represent a two-dimensional or three-dimensional surface, respectively. The two-dimensional B-spline and NURBS surfaces are defined by:

$$\boldsymbol{F}(\xi_1, \xi_2) = \sum_{(i,j) \in \mathcal{R}} \hat{B}_{(i,j)}(\xi_1, \xi_2) \boldsymbol{P}_{(i,j)} \text{ and}$$

$$\boldsymbol{F}(\xi_1, \xi_2) = \sum_{(i,j) \in \mathcal{R}} \hat{R}_{(i,j)}(\xi_1, \xi_2) \boldsymbol{P}_{(i,j)}, \text{ respectively,} \tag{6}$$

where the geometric mapping is $\boldsymbol{F}: \hat{\Omega} \to \Omega$ and $\hat{\Omega} = [0,1]^{d=2}$ is the parametric domain. The volumetric mappings are formulated in the same way. For details of B-splines and NURBS parameterization for IGA, readers may follow references [20, 25]. For simple notations, we convert the double index $(i,j)$ to a single index $k \in \mathcal{R}$. Then Eq. (6) becomes:

$$\boldsymbol{F}(\xi_1, \xi_2) = \sum_{k \in \mathcal{R}} \hat{B}_k(\xi_1, \xi_2) \boldsymbol{P}_k \text{ and}$$

$$\boldsymbol{F}(\xi_1, \xi_2) = \sum_{k \in \mathcal{R}} \hat{R}_k(\xi_1, \xi_2) \boldsymbol{P}_k, \text{ respectively.} \tag{7}$$

Based on this geometric mapping, IGA approximates the solution space in an isoparametric manner, that is, the same basis functions are used for both geometry and solution fields. Thus for the solution field, the control points are replaced with _control variables_, $c_k$, and the following solution field is given by $u^{h,IGA}(\xi_1, \xi_2) = \sum_{k \in \mathcal{R}} \hat{R}_k(\xi_1, \xi_2) c_k$ for NURBS mapping (the one with B-splines is straightforward). We assumed both $u^h$ and $c_k$ are scalar variables but they can become bold-faced to denote vector quantities: $\boldsymbol{u}^h, \boldsymbol{c}_k$. The solution field of IGA belongs to the discrete space $V_h$ in the physical domain $\Omega$:

$$V_h = \text{span}\{R_k(\mathbf{x}) \coloneqq \hat{R}_k \circ \boldsymbol{F}^{-1}(\mathbf{x}), k \in \mathcal{R}\}. \tag{8}$$

where $\boldsymbol{F}^{-1}(\mathbf{x}): \Omega \to \hat{\Omega}$ is the inverse geometric mapping and the symbol $\circ$ represents function composition.

### 2.1.2. Multi-Patch IGA

Multi-patch IGA is formulated by expanding the NURBS mapping to multiple non-overlapping and fully matching patches (B-splines follow the same procedure, but they are omitted here for brevity).

Assume the two-dimensional physical domain $\Omega$ consists of $N$ NURBS patches whose mappings are denoted by:

$$F^{(\alpha)}: \widehat{\Omega} \to \Omega^{(\alpha)}, \qquad \alpha = 1, \ldots, N \tag{9}$$

such that

$$\Omega = \bigcup_{\alpha=1}^{N} \Omega^{(\alpha)}. \tag{10}$$

The superscript $(\alpha)$ denotes functions or variables that differ by patches. Here, $\Omega^{(\alpha)}$ is a closed set. The boundary of a patch $\Omega^{(\alpha)}$ is denoted by $\Gamma^{(\alpha)}$. The $\Gamma$ represents the boundary of the whole computational domain $\Omega$ which is imposed with boundary conditions.

Each patch is connected to one or several neighboring patches. Patches $\Omega^{(\alpha)}$ and $\Omega^{(\beta)}$ ($\alpha < \beta$) are called *neighboring patches* if $\Omega^{(\alpha)} \cap \Omega^{(\beta)}$ has positive $(d-1)$-dimensional measure. Their interface is denoted by $\Gamma^{(\alpha,\beta)} = \Gamma^{(\alpha)} \cap \Gamma^{(\beta)}, \alpha < \beta$. The set $\mathcal{F}$ consists of indices of all neighboring patches, i.e., $\mathcal{F} = \{(\alpha,\beta) |\ \Omega^{(\alpha)}$ and $\Omega^{(\beta)}$ are *neighboring patches*, $\alpha < \beta\}$. Each mapping has its own NURBS basis, $\widehat{R}_k^{(\alpha)}(\xi_1, \xi_2)$. Similar to Eq. (8), the discretized solution space in the physical domain at patch $(\alpha)$ is defined as

$$V_h^{(\alpha)} = \text{span}\left\{R_k^{(\alpha)}(\mathbf{x}) := \widehat{R}_k^{(\alpha)} \circ {F^{(\alpha)}}^{-1}(\mathbf{x}), \mathbf{k} \in \mathcal{R}^{(\alpha)}\right\}, \tag{11}$$

where the patch-wise inverse mapping ${F^{(\alpha)}}^{-1}(\mathbf{x})$ is uniquely defined from the bijective assumption.

---

**Bijective Assumption**: For each patch $\Omega^{(\alpha)}$, the absolute value of the Jacobian determinant assumed to be always greater than a certain positive constant. That is, there exists $\delta > 0$ such that:

$$\left|\det\left(\nabla(F^{(\alpha)})\right)\right| > \delta, \qquad \alpha = 1, 2, \ldots, N.$$

---

In multi-patch IGA theory, the $C^0$ continuity of the solution field across patch interfaces (i.e., $G^0$ continuity – globally $C^0$) is guaranteed if the patches are fully matching. To achieve fully matching geometrical mappings, the conformity assumption is required [25, 26]:

**Conformity Assumption for IGA**: Let $\Gamma^{(\alpha,\beta)} = \Gamma^{(\alpha)} \cap \Gamma^{(\beta)}, (\alpha, \beta) \in \mathcal{F}$ be the interface between the neighboring patches $\Omega^{(\alpha)}$ and $\Omega^{(\beta)}$. Also, we define a set of indices of patch $\alpha$ whose support domain intersects the interface $\Gamma^{(\alpha,\beta)}$ as $\mathcal{B}^{(\alpha,\beta)} = \left\{\mathbf{k} \in \mathcal{R}^{(\alpha)}: \text{supp}\left(R_k^{(\alpha)}\right) \cap \Gamma^{(\alpha,\beta)} \neq \emptyset\right\}$. The two patches fully match if the two following conditions are met.

(i) The interface $\Gamma^{(\alpha,\beta)}$ is the image of a full edge of the respective parameter domains.
(ii) For each index $\mathbf{k} \in \mathcal{B}^{(\alpha,\beta)}$, there exists a unique index $\mathbf{l} \in \mathcal{B}^{(\beta,\alpha)}$, such that:

$$R_k^{(\alpha)}|_{\Gamma^{(\alpha,\beta)}} = R_l^{(\beta)}|_{\Gamma^{(\alpha,\beta)}}.$$

Given the fully matching multi-patch geometry, achieving $G^0$ continuity in the solution field is straightforward; the two matching basis functions $R_k^{(\alpha)}$ and $R_l^{(\beta)}$ should have the same control variables $c_k^{(\alpha)} = c_l^{(\beta)}$ and control weights $w_k^{(\alpha)} = w_l^{(\beta)}$. Since the control variables are the degrees of freedom (DoFs) of IGA, the imposition of $c_k^{(\alpha)} = c_l^{(\beta)}$ is simply conducted by generating a connectivity array like FEM [25].

## 2.2. Convolution Hierarchical Deep-learning Neural Network (C-HiDeNN)

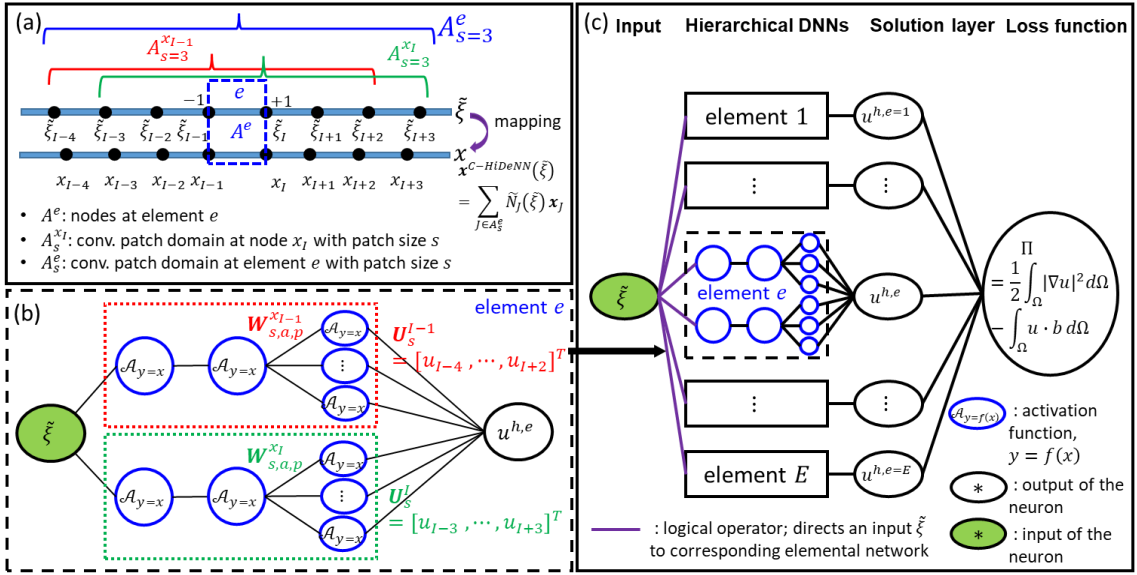

Figure 1 (a) Neural network structure of Convolution Hierarchical Deep-learning Neural Network (C-HiDeNN) for solving one-dimensional (1D) Poisson's equation. (a) Mapping between 1D parent coordinate $(\tilde{\xi})$ and physical coordinate $(x)$ with definitions on nodal and elemental convolution patches. The convolution patch size $s$ determines the number of neighboring elements employed for building *convolution patch functions*. (b) Elemental neural network for element $e$. The network is partially connected where the first two hidden layers are for linear interpolation and the third hidden layer is for nonlinear interpolation. The weights between the second and the third hidden layers are from the convolution patch functions $W_{s,a,p}^{x_I}$, while those from the third and the output layers are the corresponding nodal variables $U_s^I$. (c) The global neural network for solving 1D Poisson's equation. The purple links between the input layer and the Hierarchical DNNs denote a logical operation that directs a given input $(\tilde{\xi})$ toward a corresponding elemental network. Blue neurons with functions inside refer to the activation function. The figure is borrowed from [17] with modifications.

Convolution Hierarchical Deep-learning Neural Network (C-HiDeNN) [17-19] is a systematically pruned and partially connected neural network for higher-order interpolation without increasing the number of nodes (or DoFs). It embeds classical interpolation theories such as Lagrange polynomials/FEM [27], reproducing kernel [28, 29], and meshfree methods [30, 31] into deep-learning, the universal approximator with arbitrary choice or combination of activation functions. The neural network structure of C-HiDeNN is illustrated in Figure 1. C-HiDeNN interpolation at parent coordinate $\tilde{\xi}$ of element $e$ is given as:

$$u^{h,e}(\tilde{\xi}) = \sum_{I \in A^e} N_I(\tilde{\xi}) \sum_{K \in A_s^{x_I}} W_{s,a,p,K}^{x_I}\left(x^{h,e}(\tilde{\xi})\right) u_K = \sum_{J \in A_s^e} \tilde{N}_J(\tilde{\xi}) u_J, \qquad (12)$$

where $\tilde{\xi}$ and $x$ are the coordinates of the parent domain and the physical domain, respectively. Nodal sets $A^e, A_s^{x_I}, A_s^e$ refer to the nodes in element $e$, the convolution patch nodes centered at node $I$ with a convolution patch size $s$, and the convolution patch nodes of element $e$, respectively (illustrated in Figure 1(a) for one-dimensional case). It is worth noting here that the word "patch" refers to two different but important concepts in this paper. One is the "CAD patch," denoting B-splines or NURBS parameterized mapping $\boldsymbol{F}^{(\alpha)}$ and the corresponding mapped domains $\Omega^{(\alpha)}, \alpha = 1, 2, \ldots, N$. This will be referred to as "patch" for brevity. The other is the "convolution patch," the domain of elements that are employed for building C-HiDeNN interpolation functions, $\widetilde{N}_J(\tilde{\xi})$ in Eq. ( 12 ). To distinguish from the CAD-patch, the latter will be denoted by its full name, "convolution patch."

The *convolution patch size* $s$ determines the number of neighboring elements employed for building the third hidden layer in Figure 1(b) whose weights are called *convolution patch functions*, $\boldsymbol{W}_{s,a,p}^{x_I}(x) = \{W_{s,a,p,K}^{x_I}(x) | K \in A_s^{x_I}\}$. They are defined on a nodal convolution patch $A_s^{x_I}$, a set of nodes inside $s$ layers of elements surrounding the center node $I$ (or the corresponding physical coordinate $x_I$). Since each element has $n(A^e)$ nodes, the support domain of element $e$ (i.e., $A_s^e$) is the union of all nodal convolution patches: $A_s^e = \cup_{I \in A^e} A_s^{x_I}$. For example, in Figure 2(c), the elemental convolution patch of the blue element $e$ with $s = 1$ designates the green elements, one layer of elements surrounding the blue element. The same rule applies to $s = 2$ in Figure 2(d). If an element is near the domain boundary, the elemental convolution patch $A_s^e$ narrows because no more elements exist beyond the boundary. See Figure 2(e-f).

The function $N_I(\tilde{\xi})$ in Eq. ( 12 ) is the linear shape function at node $I$, and $W_{s,a,p,K}^{x_I}(x)$ is the convolution patch function for node $K$ inside the nodal convolution patch $A_s^{x_I}$. Note that the linear interpolant can be any higher order elements [1], but will be kept linear throughout this paper as the nonlinear part can be more systematically handled by the convolution patch functions. They are constructed with physical coordinates based on the meshfree interpolation theories [17], although any other combinations with neural network activation functions and polynomial functions can also be used [32]. A detailed procedure of constructing the convolution patch functions $\boldsymbol{W}_{s,a,p}^{x_I}(x)$ is provided in Appendix A. They are parameterized with the *dilation parameter* $a$ that determines the influence domain of the radial basis function and the *reproducing polynomial order* $p$, as well as the *convolution patch size* $s$. Since they are built with physical nodes, the parent coordinate $\tilde{\xi}$ is mapped to the physical coordinate $x$ via the FEM mapping, $x^{h,e}(\tilde{\xi}) = \sum_{I \in A^e} N_I(\tilde{\xi}) x_I$, then fed into the function $W_{s,a,p,K}^{x_I}\left(x^{h,e}(\tilde{\xi})\right)$. Finally, the two functions $N_I(\tilde{\xi})$ and $W_{s,a,p,K}^{x_I}(x)$ are combined to form C-HiDeNN interpolants $\widetilde{N}_J(\tilde{\xi})$ at node $J \in A_s^e$.

As illustrated in Figure 2, once the shape functions are constructed, C-HiDeNN follows the standard procedure of interpolation: nodal shape function times nodal variable, followed by a summation over nodes. Compared to standard FEM (linear or higher order, see Figure 2(a-b)) that uses only the nodes inside an element to approximate a solution field of that element, C-HiDeNN uses nodes of neighboring elements, as illustrated in Figure 2(c-f).

C-HiDeNN can build higher order shape functions using only linear elements (3- or 4-node elements in 2D and 4- or 8-node elements in 3D), achieving faster convergence rates than FEM without increasing DoFs. Instead, C-HiDeNN broadens the bandwidth of the global stiffness matrix due to the extended nodal connectivity. The computational overhead coming from the computation of convolution patch functions is resolved by graphics processing unit (GPU) programming because this operation is parallelizable. Detailed numerical studies can be found in [17-19].

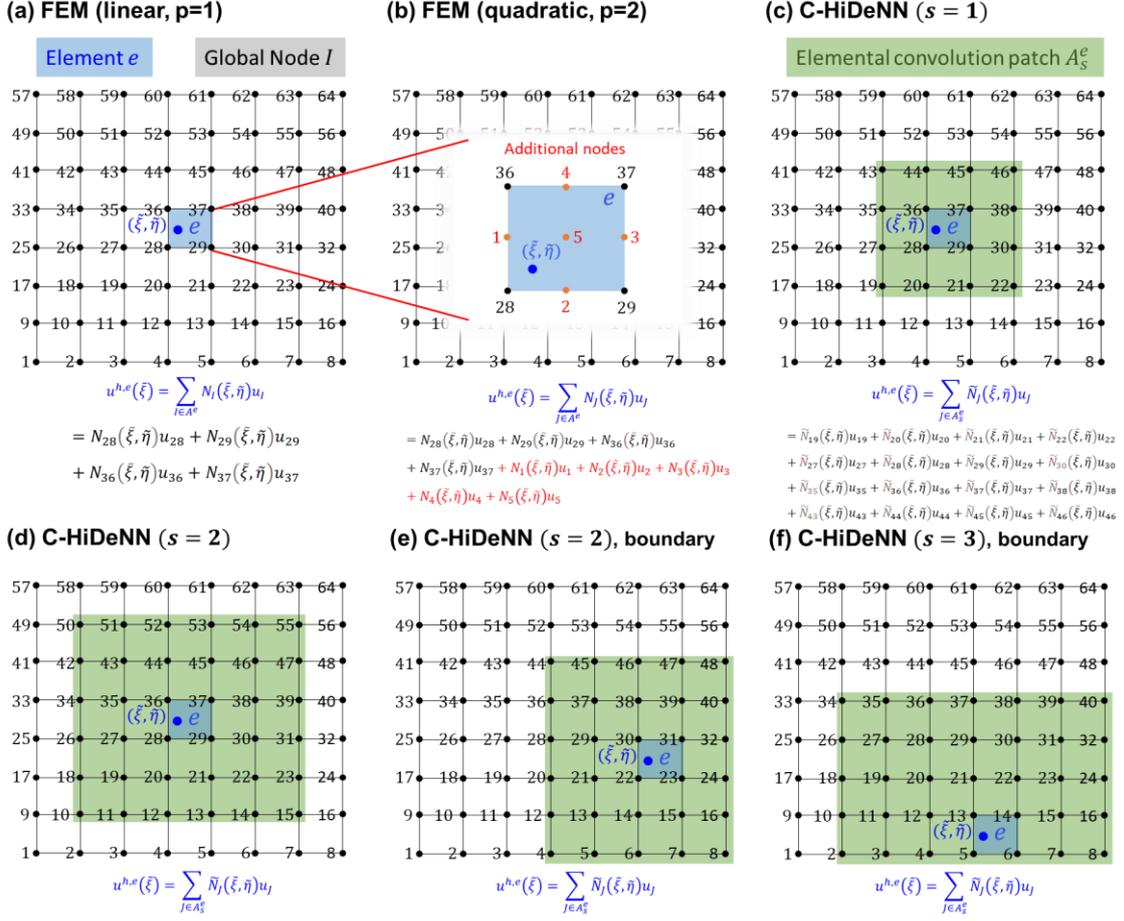

Figure 2 Schematic illustration of nodal interpolation of (a-b) FEM and (c-f) C-HiDeNN. FEM and C-HiDeNN interpolation functions are $N_I(\tilde{\xi}, \tilde{\eta})$ and $\widetilde{N}_J(\tilde{\xi}, \tilde{\eta})$, respectively.

### 2.3. Single-Patch Isogeometric C-HiDeNN (in short, C-IGA)

Because C-HiDeNN is built upon linear FEM mesh, it cannot reproduce an exact geometry like what IGA does. However, if the convolution patch function is enforced to reproduce NURBS (or B-spline) basis function instead of polynomials, one can create a C-HiDeNN geometrical mapping that reproduces the exact geometry like what IGA does. This idea was first introduced in [2] and named as Isogeometric C-HiDeNN or in short, C-IGA. The authors proved the equivalence of the geometric mappings of IGA and C-IGA, which is (when $d = 2 \text{ or } 3$),

$$\boldsymbol{F}(\boldsymbol{\xi}) = \sum_{k \in \mathcal{R}} \hat{R}_k(\boldsymbol{\xi}) \boldsymbol{P}_k = \sum_{J \in A_s^e} \widetilde{N}_J(\boldsymbol{\xi}) \boldsymbol{x}_J, \text{where } \boldsymbol{x}_J = \boldsymbol{F}(\boldsymbol{\xi}_J). \qquad (13)$$

where the IGA geometric mapping follows Eq. ( 7 ). Note that C-IGA convolution patch functions $W_{s,a,p}^{\xi_I}(\boldsymbol{\xi})$ are created in the parametric domain $\boldsymbol{\xi} \in \widehat{\Omega}$, whereas those of C-HiDeNN $W_{s,a,p}^{x_I}(\boldsymbol{x})$ are defined in the physical domain $x \in \Omega$. For C-IGA interpolation, we need to discretize the problem domain similar to FEM and the discrete nodes are defined as $x_J$ and $\xi_J$ for the physical and parametric domains, respectively. These nodes follow the IGA forward mapping: $x_J = F(\xi_J)$. The key difference between IGA and C-IGA arises here. Considering both methods adopt the isoparametric formulation, IGA solution field $u^{h,IGA}(\xi)$ interpolates (or regresses) control variables $c_k$ attached to the control points $P_k$, whereas C-IGA solution field $u^{h,C-IGA}(\xi)$ interpolates nodal variables $u_J$ attached to $x_J$, as illustrated in Figure 3. Note that the solution field $u^h$, control variables $c_k$, and nodal variables $u_J$ become bold-faced ($\boldsymbol{u^h}, \boldsymbol{c_k}, \boldsymbol{u_J}$) when there are multiple output variables (e.g., x-,y-,z-displacements at node $J$).

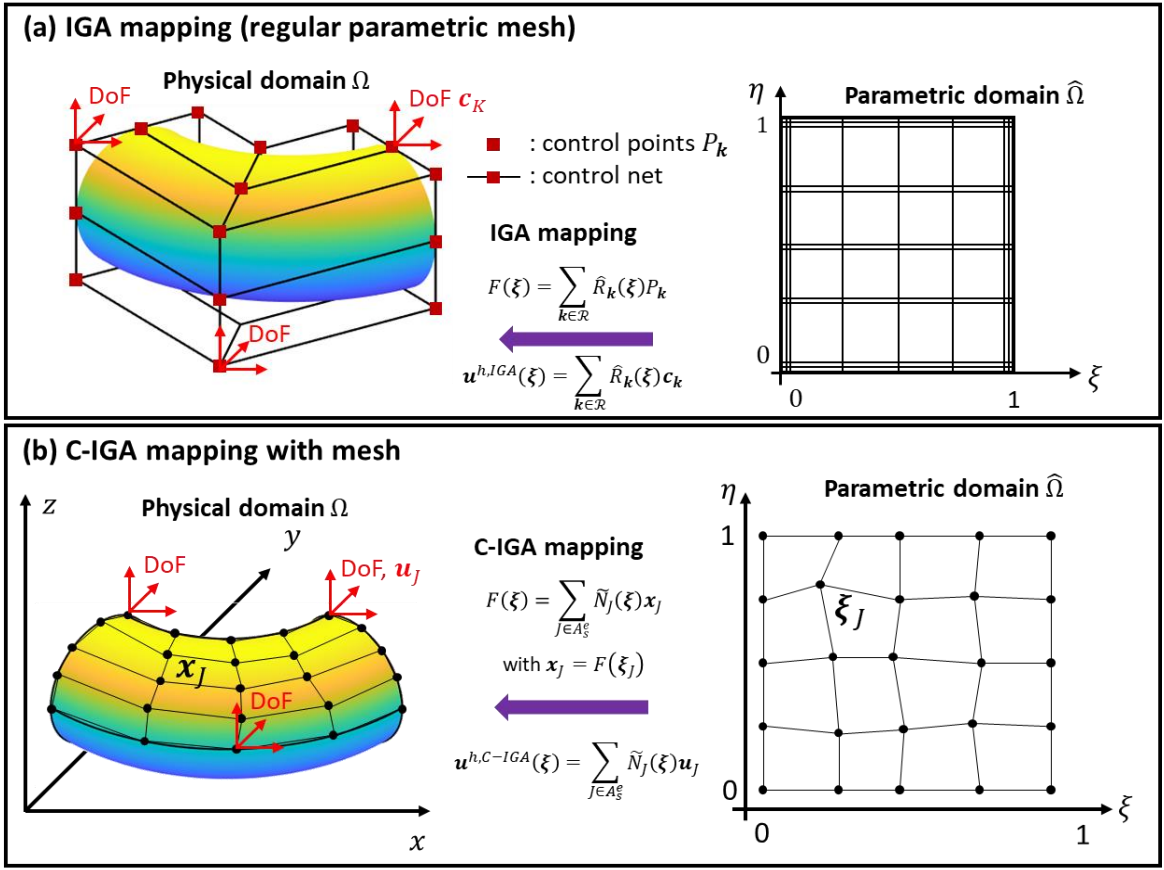

Figure 3 Schematic illustration of the geometric mappings of IGA (a) and C-IGA (b). The figure is borrowed from [2] with modifications.

The C-IGA mesh can be created by two opposite ways: 1) generate parametric mesh $(\xi_J)$ and convert it to the physical mesh $(x_J)$; 2) the opposite, from the physical mesh to the parametric mesh. The former is easier to program because the forward geometric mapping $F: \xi \rightarrow x$ is provided from a CAD model. However, the opposite (the latter) is more practical and promising if one wants to refine a mesh in the physical domain like FEM. In this case, the inverse mapping $F^{-1}: x \rightarrow \xi$ should be implemented, which is discussed further in Section 4, Algorithm 1. Under the isoparametric formulation, the C-IGA solution field at element $e$ is given by:

$$u^{h,e} = \sum_{J \in A_s^e} \widetilde{N}_J(\xi) u_J \qquad (14)$$

where $u_J$ is the nodal variable corresponding to the nodal coordinate $x_J$ and $\xi_J$.

C-IGA inherits advantages of C-HiDeNN such as versatile adaptivity (r-, h-, s-, a-, p-), Kronecker delta property, and higher-order approximation without DoFs increase, while achieving the exact geometry of a CAD model. Table 2 summarizes the differences between higher order FEM, IGA, physics-informed neural network (PINN), C-HiDeNN, and C-IGA. The key difference to be highlighted is that C-IGA interpolates nodal variables $(u_J)$ attached to the physical geometry, while IGA cannot interpolate control variables $(c_k)$ in general and they are not attached to the physical geometry for most cases (see Figure 3(a)). This requires additional treatment for imposing boundary conditions in IGA, whereas C-IGA can readily follow the standard procedure of FEM. Furthermore, C-IGA can be formulated on arbitrary or unstructured meshes in the physical domain $\Omega$ (see Figure 3(b)), rendering versatile mesh adaptivity like standard FEM. On the other hand, the mesh refinement of IGA is mostly conducted in the parametric domain $\widehat{\Omega}$ with a global refinement for NURBS [20] or a local refinement for T-splines [21]. As illustrated in Figure 4, an unstructured mesh consisting of quadrilateral and triangular elements in 2D can be generated in the physical domain (a) and then inverse-mapped into the parametric domain (b). C-IGA can solve the problem with this unstructured mesh without sacrificing accuracy compared to C-IGA or IGA with regular mesh (c). Detailed numerical studies for single-patch C-IGA can be found in [2].

Another highlight drawn in Table 2 is the comparison of PINN with the others. In PINN, the approximate solution field is constructed via neural networks with adaptable number of neurons, hidden layers, and activation functions [5]. Nevertheless, there are no such concepts like discrete mesh, nodal variables, interpolation functions, and reproducing properties, making it difficult to interpret and conduct convergence study of the numerical solution. On the other hand, C-HiDeNN/C-IGA are interpretable neural networks resembling FEM and IGA because the networks are constructed based on the interpolation theory as illustrated in Figure 1 and in Eq. ( 12 ). They can achieve higher order interpolations and the convergence of error estimates can be proven mathematically [2, 17, 18].

Table 2 Comparison between FEM, IGA, physics-informed neural network (PINN), C-HiDeNN, and C-IGA.

|  | FEM | IGA | PINN | C-HiDeNN | C-IGA |
|---|---|---|---|---|---|
| DoFs attached to | nodes | control variables | weights, biases | nodes | nodes |
| Interpolation (Kronecker delta) | Yes | In general, No | No | Yes | Yes |
| Reproducing property | polynomials up to order $p$ | exact geometry | No | any nonlinear activation | exact geometry with reproducing order $p$ |
| Adaptivity | $r, h, p$ | $h, k$ | neurons, layers, | $r, h, p, s, a$ | $r, h, p, s, a$ |

activations

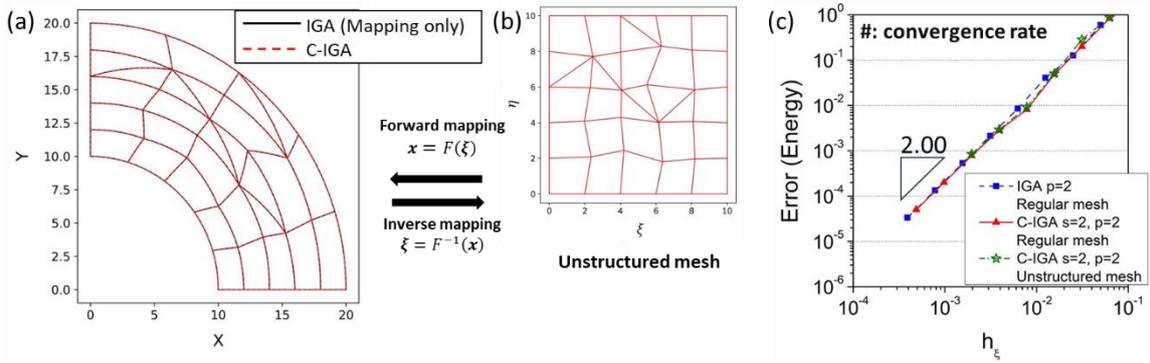

Figure 4 Schematic illustration of unstructured mesh created to solve 2D Poisson's equation using C-IGA. Physical domain (a), parametric domain (b), and convergence plot with convergence rate (c). IGA can only do the mapping; it cannot solve the problem with this mesh whereas C-IGA can. Detailed problem definitions can be found in [2].

## 3. Multi-Patch C-IGA Formulation without Compatibility Conditions

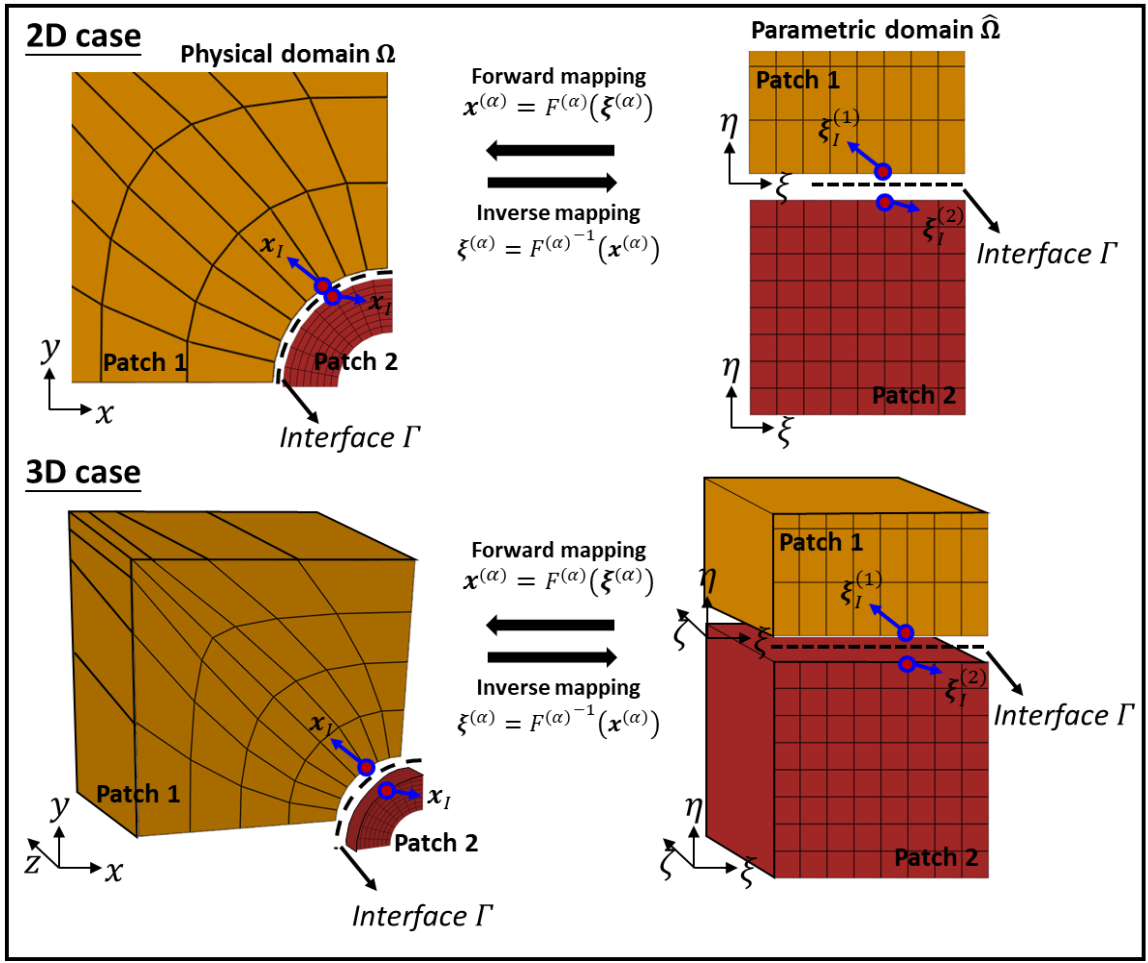

Figure 5 Two-dimensional (top) and three-dimensional cases (bottom) for a two-CAD-patch system. Physical space (left) and parametric space (right). Yellow box represents patch $\alpha$ ($\Omega^{(\alpha)}$), and red box represents patch $\beta$ ($\Omega^{(\beta)}$).

Within a single CAD patch, C-IGA achieves $C^0$ continuity at the element boundary and higher order of continuity inside each element [2]. On the other hand, for a multi-patch system stated in Section 2.1.2, C-IGA interpolation may suffer from compatibility issues at patch interfaces, similar to IGA without the **Conformity Assumption** (see Section 2.1.2). In this section, we discuss the effect of this non-conforming patch interface on the convergence rate of a general interpolation.

First, we define a general interpolation operator $\mathcal{J}$. The interpolation of a continuous function in a Sobolev space $H^{p+1}(\Omega)$, $u(x) \in H^{p+1}(\Omega)$ ($p \geq 0$) reads:

$$u^h(x) = \mathcal{J}u(x), \qquad (15)$$

where the interpolated function $u^h(x)$ is constructed with multiple non-overlapping patches, each of which is parameterized differently. Thus, $u^h(x)$ is patch-wise continuous (i.e., $u^{h,(\alpha)} \in H^1(\Omega^{(\alpha)})$) and can be discontinuous at patch interfaces due to incompatible parameterizations. We define the following *broken Sobolev space* [33]:

$$H_b^k(\Omega) = \{v \in L^2(\Omega): v^{(\alpha)} \in H^k(\Omega^{(\alpha)}), \alpha = 1,2,\dots,N\}. \qquad (16)$$

where $v^{(\alpha)}$ is a function confined in $\Omega^{(\alpha)}$. Obviously, $u^h \in H_b^{p+1}(\Omega)$ and $H^k(\Omega)$ is a subset of $H_b^k(\Omega)$. With this setup, we investigate the convergence of the interpolation of a continuous function in Section 3.1 and the solution of a partial differential equation (PDE) in Section 3.2.

**3.1 Interpolation Estimate on a Multi-Patch System without Compatibility Conditions**

We define a norm in $H_b^k(\Omega)$ as:

$$\|v\|_{H_b^k(\Omega)} = \left(\left(\|v^{(1)}\|_{H^k(\Omega^{(1)})}\right)^2 + \cdots + \left(\|v^{(N)}\|_{H^k(\Omega^{(N)})}\right)^2\right)^{1/2}. \qquad (17)$$

If the interpolation $u^h = \mathcal{J}u$ can reach the $p$-th order interpolation estimate in each patch $\Omega^{(\alpha)}$, the interpolation of a continuous function $u$ throughout the whole domain $\Omega$ can also reach the $p$-th order convergence rate. See Lemma 1.

**Lemma 1**:
For $u(x) \in H^{p+1}(\Omega)$, assume that the $p$-th order interpolation estimate for $u^h = \mathcal{J}u$ holds in each single patch $\Omega^{(\alpha)}, \alpha = 1,2,\dots,N$, i.e.,

$$\|u^{(\alpha)} - u^{h,(\alpha)}\|_{H^1(\Omega^{(\alpha)})} \leq C_{interp}^{(\alpha)} h^p \|u\|_{H^{p+1}(\Omega^{(\alpha)})}, \alpha = 1,2,\dots,N, \qquad (18)$$

where $C_{interp}^{(\alpha)}$ is a constant independent of the average element size $h$. Then the following interpolation estimate in the whole domain $\Omega$ holds:

$$\|u - u^h\|_{H_b^1(\Omega)} \leq C_{interp} h^p \|u\|_{H^{p+1}(\Omega)} \qquad (19)$$

with $C_{interp} = \max\left(C_{interp}^{(1)}, C_{interp}^{(2)}, \dots, C_{interp}^{(N)}\right)$.

The proof is provided in the Appendix, and numerical illustrations are provided in Sections 5.1 and 5.2.

**3.2 Convergence Rates for Solving Elliptic PDEs without Compatibility Conditions**

In contrast to the interpolation estimate of a continuous function discussed in the previous section, the convergence rate of a solution to a PDE is suboptimal due to the compatibility issue at patch interfaces. Here, we take an elliptic PDE to illustrate the effect of compatibility.

Consider the following elliptic PDE:
$$-\Delta u + u = f \quad (20)$$

with $f$ being a source term. The boundary $\Gamma$ of the entire domain $\Omega$ consists of Dirichlet boundary $\Gamma^D$ assigned with $u|_{\Gamma^D} = u_D$, and Neumann boundary $\Gamma^N$ assigned with $\left.\frac{\partial u}{\partial n}\right|_{\Gamma^N} = g$. Note that $\Gamma^D \neq \emptyset$ and $\Gamma^D, \Gamma^N$ do not overlap.

The solution $u \in H^{p+1}(\Omega)$ of Eq. (20) satisfies the following strong form PDE for any test function $w \in H_b^1(\Omega)$,

$$\sum_{\alpha=1}^{N} \int_{\Omega^{(\alpha)}} (-\Delta u + u - f) w \, dV = 0, \quad \forall w \in H_b^1(\Omega). \quad (21)$$

By integration by parts, Eq. (21) becomes:

$$a(w, u) - \sum_{\Gamma^{(\alpha,\beta)} \in \mathcal{F}} \left([w]_{\Gamma^{(\alpha,\beta)}}, \langle \frac{\partial u}{\partial n^{(\alpha)}} \rangle_{\Gamma^{(\alpha,\beta)}}\right)_{\Gamma^{(\alpha,\beta)}} = (w, f) + (w, g)_{\Gamma^N}, \quad (22)$$

$$\forall w \in H_b^k(\Omega).$$

The $a(\cdot, \cdot)$ is a bilinear form given by:

$$a(w, v) = \sum_{\alpha=1}^{N} \int_{\Omega^{(\alpha)}} (\nabla w \cdot \nabla v + w \cdot v) dV(x), \quad w, v \in H_b^k(\Omega). \quad (23)$$

Note that $\left(a(w, w)_{\Omega^{(\alpha)}}\right)^{\frac{1}{2}}$ is precisely the $H_b^1$ norm of $w$ over the whole domain $\Omega$, that is, $\|w\|_{H_b^1(\Omega)}$. Here, for convenience, we introduce the following notations for $u \in H_b^1(\Omega)$:

- $\langle u \rangle_{\Gamma^{(\alpha,\beta)}} = \frac{1}{2}\left(u^{(\alpha)}|_{\Gamma^{(\alpha,\beta)}} + u^{(\beta)}|_{\Gamma^{(\alpha,\beta)}}\right)$ – mean value of the traces of $u$ at the interface $\Gamma^{(\alpha,\beta)}$.

- $[u]_{\Gamma^{(\alpha,\beta)}} = u^{(\alpha)}|_{\Gamma^{(\alpha,\beta)}} - u^{(\beta)}|_{\Gamma^{(\alpha,\beta)}}$ – jump of $u$ at the interface $\Gamma^{(\alpha,\beta)}$.

Remark that for other elliptic PDEs, one may find the corresponding equivalent norm $\|\cdot\|$ in $H_b^k(\Omega)$ like the one in the single patch [27], i.e.,

$$c_1 \|w\| \leq \left(a(w, w)\right)^{\frac{1}{2}} \leq c_2 \|w\|, \quad (24)$$

where $c_1$ and $c_2$ are positive constants. Then the error in the following analysis will be measured by this equivalent norm $\|\cdot\|$.

The weak form in a multi-patch domain $\Omega$ reads:

$$a(w^h, u^h) + \rho \sum_{\Gamma^{(\alpha,\beta)} \in \mathcal{F}} \left([w^h]_{\Gamma^{(\alpha,\beta)}}, [u^h]_{\Gamma^{(\alpha,\beta)}}\right)_{\Gamma^{(\alpha,\beta)}} = (w^h, f) + (w^h, g)_{\Gamma^N}, \tag{25}$$

$$\forall w^h \in \mathcal{V}^h,$$

where $u^h \in \mathcal{S}^h$ is the trial function in $H^{p+1}(\Omega)$ with the Dirichlet boundary conditions on $\Gamma^D$, i.e.,

$$\mathcal{S}^h = \{v^h | v^h = \mathcal{I}v, v \in H^{p+1}(\Omega) \text{ and } v|_{\Gamma^D} = u_D\}, \tag{26}$$

and $\mathcal{V}^h$ is the set of test functions in $H^{p+1}(\Omega)$ with zero boundary conditions on $\Gamma^D$, i.e.,

$$\mathcal{V}^h = \{v^h | v^h = \mathcal{I}v, v \in H^{p+1}(\Omega) \text{ and } v|_{\Gamma^D} = 0\}. \tag{27}$$

The term $\rho \sum_{\Gamma^{(\alpha,\beta)} \in \mathcal{F}} \left([w^h]_{\Gamma^{(\alpha,\beta)}}, [u^h]_{\Gamma^{(\alpha,\beta)}}\right)_{\Gamma^{(\alpha,\beta)}}$ is the penalty term to guarantee the uniqueness of the solution, where $\rho$ is a selected positive penalty coefficient.

The weak form of a single patch system will lead to $p$-th order convergence rate if the $p$-th order interpolation estimate is achieved. However, when solving a multi-patch system, the compatibility issue will limit the convergence according to Theorem 1.

---

**Theorem 1**:

Let $u^h \in \mathcal{S}^h$ be the trial function of the weak form:

$$a(w^h, u^h) + \rho \sum_{\Gamma^{(\alpha,\beta)} \in \mathcal{F}} \left([w^h]_{\Gamma^{(\alpha,\beta)}}, [u^h]_{\Gamma^{(\alpha,\beta)}}\right)_{\Gamma^{(\alpha,\beta)}} = (w^h, f) + (w^h, g)_{\Gamma^N}, \tag{28}$$

$$\forall w^h \in \mathcal{V}^h,$$

where $\rho > 0$ is chosen large enough to guarantee the uniqueness of the solution.
Assume that the $p$-th order interpolation estimate for any interpolation $v^h \in \mathcal{S}^h$ holds. If an estimate of the discontinuity at the interface $\Gamma^{(i,j)}$ satisfies:

$$\left\| [w^h]_{\Gamma^{(\alpha,\beta)}} \right\|_{L^2(\Gamma^{(\alpha,\beta)})} \leq C_{gap} \left(h^{\Gamma^{(\alpha,\beta)}}\right)^q \left( \left\|w^{h,(\alpha)}\right\|_{L^2(\Gamma^{(\alpha,\beta)})} + \left\|w^{h,(\beta)}\right\|_{L^2(\Gamma^{(\alpha,\beta)})} \right), \tag{29}$$

$$\forall w^h \in \mathcal{V}^h,$$

where $h^{\Gamma^{(\alpha,\beta)}}$ is the element size around the interface $\Gamma^{(\alpha,\beta)}$, $q \in \mathbb{N}$ and $C_{gap} > 0$ are constants, then the following error estimate holds:

$$\|u^h - u\|_{H_b^1(\Omega)} \leq \sum_{\Gamma^{(\alpha,\beta)} \in \mathcal{F}} C_1^{\Gamma^{(\alpha,\beta)}} C_{gap} \left(h^{\Gamma^{(\alpha,\beta)}}\right)^q \left\| \langle \frac{\partial u}{\partial \mathbf{n}^{(\alpha)}} \rangle_{\Gamma^{(\alpha,\beta)}} \right\|_{L^2(\Gamma^{(\alpha,\beta)})} + C_2 h^p \|u\|_{H^{p+1}(\Omega)}, \tag{30}$$

with $u$ being the solution to Eq. (20) under the Dirichlet and Neumann boundary conditions.

---

**Corollary:** If the interpolations always achieve compatibility at the interfaces, i.e., the trial function $u^h \in \mathcal{S}^h$ and the test function $w^h \in \mathcal{V}^h$ are both continuous at the interface: $[u^h]_{\Gamma^{(\alpha,\beta)}} = 0, [w^h]_{\Gamma^{(\alpha,\beta)}} = 0$, then the coefficient $C_{gap}$ in Eq. (29) becomes zero. Accordingly, the first term on the right-hand side of inequality (30) vanishes, and the following error estimate holds:

$$\|u^h - u\|_{H_b^1(\Omega)} \leq C_2 h^p \|u^{Ext}\|_{H^{p+1}(\Omega)}, \tag{31}$$

that leads to the optimal convergence of the error estimate.

**Theorem 1** indicates that the convergence rate is limited by the gap of interpolations at the interface. In general, $q$ is smaller than $p$. For C-IGA interpolations, the numerical tests provided in Section 5 show that $q = 1$. Hence, the final convergence rate without compatibility is merely first order even though the reproducing order of interpolation for each patch is greater than one. Thus, one must satisfy the *compatibility conditions* for multi-patch systems to ensure $C^0$ continuity globally (i.e., $G^0$ continuity), and achieve the theoretical convergence rate.

## 4. Multi-Patch C-IGA Formulation with Compatibility Conditions

The single-patch C-IGA interpolation at patch $\Omega^{(\alpha)}$ is defined as

$$u^{h,(\alpha)} = \sum_I N_I^{(\alpha)}(\boldsymbol{\xi}^{(\alpha)}) \sum_{K \in A_s^{\xi_I^{(\alpha)}}} W_{s,a,p,K}^{\xi_I^{(\alpha)}}(\boldsymbol{\xi}^{(\alpha)}) u_K = \sum_{J \in A_s^e} \widetilde{N}_J^{(\alpha)}(\boldsymbol{\xi}^{(\alpha)}) u_J \tag{32}$$

with an isoparametric mapping

$$\boldsymbol{x}(\boldsymbol{\xi}^{(\alpha)}) = \sum_I N_I^{(\alpha)}(\boldsymbol{\xi}^{(\alpha)}) \sum_{K \in A_s^{\xi_I^{(\alpha)}}} W_{s,a,p,K}^{\xi_I^{(\alpha)}}(\boldsymbol{\xi}^{(\alpha)}) \boldsymbol{x}_J = \sum_{J \in A_s^e} \widetilde{N}_J^{(\alpha)}(\boldsymbol{\xi}^{(\alpha)}) \boldsymbol{x}_J. \tag{33}$$

The notation of involved variables is listed in Table 1. According to the single-patch C-IGA theory, the isoparametric mapping in Eq. (33) and the IGA forward mapping in Eq. (9) are equivalent. Moreover, the interpolation in Eq. (32) has two properties:

*(1)* $p$-th order reproducing property, and

*(2)* Kronecker delta property.

The former guarantees the consistency and the $p$-th order interpolation estimate in each patch. The latter imposes interpolating features to the C-IGA shape functions, which is a crucial requirement of the development of the compatibility conditions. The latter imposes interpolating features to C-IGA shape functions, which is a crucial requirement for the compatibility conditions.

**Theorem 1** indicates that the convergence rate of an elliptic PDE solution in a multi-patch system is suboptimal due to the incompatibility at the interface. To resolve this, it is required to have a $G^0$ *compatibility* defined as follows [23].

---
**$G^0$ compatibility** : The solution field achieves global $C^0$ continuity over the entire multi-patch physical domain $\Omega$, having $C^0$ continuity across patch interfaces and up to $C^N$ ($N = \min(m,n)$) continuity within each patch, where FEM shape functions are $C^n$ continuous and convolution patch functions are $C^m$ continuous [2].

---

This is realized by two steps:

*(1)* **nodal compatibility condition** to guarantee the uniqueness of solution. It replaces the penalty term in Eq. (25).

*(2)* $G^0$ **compatibility conditions** for C-IGA shape functions to achieve point-wise continuous interpolations ($C^0$ continuity) along the interface.

### 4.1 Nodal Compatibility Condition

The *nodal compatibility* gives constraints on the mesh construction. Due to the Kronecker delta property of C-IGA interpolations, if the two neighboring patches share the same physical nodes, the nodal compatibility is achieved. That is, $u^{(\alpha)}(x_I) = u^{(\beta)}(x_I)$, $\forall x_I \in \Gamma^{(\alpha,\beta)}, \Gamma^{(\alpha,\beta)} \in \mathcal{F}$, where the subscript $I$ being the nodal index along the interface. Thus, the uniqueness of the solution is guaranteed, and the penalty term can be eliminated in Eq. ( 25 ).

Achieving nodal compatibility is straightforward. First, construct a conforming physical mesh in the whole domain $\Omega$ (no hanging nodes at the interface) using a standard FEM meshing tool. Then find corresponding parametric mesh for each patch using the inverse mapping. We devised an inverse mapping algorithm based on the neural network, which will be described below.

The inverse mapping of a multi-patch system is defined as:

$$F^{(\alpha)^{-1}}: x^{(\alpha)} \to \xi^{(\alpha)}, \qquad x \in \Omega^{(\alpha)}, \alpha = 1,2,\ldots,N. \qquad (34)$$

According to the implicit function theorem, if the Jacobian determinant of a forward mapping is non-zero in the given computational domain, i.e., $\frac{\partial x^{(\alpha)}}{\partial \xi^{(\alpha)}} \neq 0$, the inverse mapping $F^{(\alpha)^{-1}}$ exists. This property is mentioned in the previous **Bijective Assumption**.

In the mesh generation step of C-IGA, we need to find out the parametric nodal coordinates $\xi_I^{(\alpha)}$ that correspond to the given physical coordinates $x_I^{(\alpha)}$ (here, the nodal index $I$ refers to the nodes in the entire patch $\Omega^{(\alpha)}$, not only in the patch interface). One possible approach is to use a feedforward neural network (FFNN) to train the inverse IGA mapping. This is a regression problem to be solved by minimizing the loss function below:

$$Loss = \frac{1}{N_L} \sum_{L=1}^{N_L} \left\| F_{NN}^{(\alpha)^{-1}}\left(x_L^{(\alpha)}\right) - \xi_L^{(\alpha)} \right\|^2 \qquad (35)$$

where the subscript $L$ denotes the index of sampling points from the IGA forward mapping $x_L^{(\alpha)} = F^{(\alpha)}\left(\xi_L^{(\alpha)}\right)$, and $F_{NN}^{(\alpha)^{-1}}$ is the FFNN inverse mapping. Once the neural network is trained successfully, we use the trained inverse mapping to generate parametric mesh $\xi_I^{(\alpha)}$ according to the given physical mesh $x_I^{(\alpha)}$. Detailed procedures of training a neural network for the inverse IGA mapping are given in Algorithm 1.

> **Algorithm 1**: Training procedure of IGA inverse mapping.
>
> - **Given:** physical mesh, $x_I^{(\alpha)}$, and IGA forward mapping, $x^{(\alpha)} = F^{(\alpha)}(\xi^{(\alpha)})$.
> - **Goal:** find parametric mesh $\xi_I^{(\alpha)}$ that satisfies the IGA forward mapping.
> 1. Generate training data. Generate $\left(\xi_L^{(\alpha)}, x_L^{(\alpha)}\right)$ pairs using the given IGA forward mapping. The training data can be different from the physical mesh.
> 2. Train a neural network for the inverse mapping, $\xi^{(\alpha)} = F_{NN}^{(\alpha)^{-1}}(x^{(\alpha)})$.
>    - A standard FFNN with sigmoid activation function and mean-squared error (MSE) loss function is a fair starting point.
> 3. For the given physical mesh, compute parametric mesh from the trained inverse mapping:
>    $$\xi_{I,NN}^{(\alpha)} = F_{NN}^{(\alpha)^{-1}}\left(x_I^{(\alpha)}\right).$$
> 4. Post processing:
>    - If $\left\|F^{(\alpha)}\left(\xi_{I,NN}^{(\alpha)}\right) - x_I^{(\alpha)}\right\| < \epsilon$ (*tolerance*), use $\xi_{I,NN}^{(\alpha)}$ as $\xi_I^{(\alpha)}$.
>    - If not, conduct a postprocessing. Solve an optimization problem with $\xi_{I,NN}^{(\alpha)}$ being the initial condition: $\xi_I^{(\alpha)} = arg\min_{\xi^{(\alpha)}}\left\|F^{(\alpha)}(\xi^{(\alpha)}) - x_I^{(\alpha)}\right\|$. Assuming that a unique inverse mapping exists, the loss of the optimization will converge to zero (with machine error).

Although C-IGA formulation with nodal compatibility is straightforward to implement, it does not guarantee $G^0$ continuity because the C-IGA shape functions of patch $(\alpha)$ are constructed only with nodes inside that patch. Thus, additional constraints should be imposed to achieve $G^0$ continuity.

### 4.2 $G^0$ Compatibility Condition

The $G^0$ *compatibility* requires further constraints on the C-IGA shape functions in addition to the nodal compatibility. The objective is to achieve point-wise continuity across the interfaces of every pair of neighboring patches $\Omega^{(\alpha)}$ and $\Omega^{(\beta)}$, that is,

$$u^{h,(\alpha)}(x) = u^{h,(\beta)}(x), \forall\, x \in \Gamma^{(\alpha,\beta)}. \tag{36}$$

In Theorem 2 below, we assert that the three compatibility conditions on the FEM shape functions, and the convolution patch functions whose support domains intersect with the interface are sufficient to establish the $G^0$ compatibility.

For ease of presentation, we make some notations. Assume that two neighboring patches $\Omega^{(\alpha)}, \Omega^{(\beta)}$ satisfy the *nodal compatibility*, that is, share the same physical nodes along the interface $\Gamma^{(\alpha,\beta)}$. Consider C-IGA interpolations in $\Omega^{(\alpha)}$ and $\Omega^{(\beta)}$ use the same parameters $s, a, p$. For a shared physical node $x_I \in \Gamma^{(\alpha,\beta)}$, we denote the corresponding parametric nodes in $\widehat{\Omega}^{(\alpha)}$ and $\widehat{\Omega}^{(\beta)}$ as $\xi_I^{(\alpha)}, \xi_I^{(\beta)}$, respectively. Let the nodal convolution patches for node $I$ be $A_s^{\xi_I^{(\alpha)}}$, $A_s^{\xi_I^{(\beta)}}$, respectively. The index set of shared nodes in the convolution patches $A_s^{\xi_I^{(\alpha)}}$, $A_s^{\xi_I^{(\beta)}}$ is denoted by $\mathcal{P}_s^{x_I} = A_s^{\xi_I^{(\alpha)}} \cap A_s^{\xi_I^{(\beta)}}$. The index set of internal nodes in $A_s^{\xi_I^{(\alpha)}}$ is $Q_s^{\xi_I^{(\alpha)}} = A_s^{\xi_I^{(\alpha)}} \setminus \mathcal{P}_s^{x_I}$, and that in $A_s^{\xi_I^{(\beta)}}$ is $Q_s^{\xi_I^{(\beta)}} = A_s^{\xi_I^{(\beta)}} \setminus \mathcal{P}_s^{x_I}$.

---

**Theorem 2**: Let two neighboring patches $\Omega^{(\alpha)}, \Omega^{(\beta)}$ **satisfy** the *nodal compatibility condition*, and the C-IGA interpolations on them use the same parameters $s, a, p$. For each shared node $I$ at the interface $\Gamma^{(\alpha,\beta)}$, if the following compatibility conditions for C-IGA interpolations hold:

**(C1)** $N_I^{(\alpha)}(x)|_{\Gamma^{(\alpha,\beta)}} = N_I^{(\beta)}(x)|_{\Gamma^{(\alpha,\beta)}}, \forall\, x \in \Gamma^{(\alpha,\beta)}$ (FEM shape functions match at the interface)

**(C2)** For $K \in \mathcal{P}_s^{x_I}$, $W_{s,a,p,K}^{\xi_I^{(\alpha)}}(x)|_{\Gamma^{(\alpha,\beta)}} = W_{s,a,p,K}^{\xi_I^{(\beta)}}(x)|_{\Gamma^{(\alpha,\beta)}}, \forall\, x \in \Gamma^{(\alpha,\beta)}$

(convolution patch functions for the shared nodes match at the interface)

**(C3)** For $K_1 \in Q_s^{\xi_I^{(\alpha)}}, K_2 \in Q_s^{\xi_I^{(\beta)}}$, $W_{s,a,p,K_1}^{\xi_I^{(\alpha)}}(x)|_{\Gamma^{(\alpha,\beta)}} = 0, W_{s,a,p,K_2}^{\xi_I^{(\beta)}}(x)|_{\Gamma^{(\alpha,\beta)}} = 0, \forall\, x \in \Gamma^{(\alpha,\beta)}$

(convolution patch functions for the internal nodes vanish at the interface),

then the C-IGA interpolation reaches $\underline{G^0\ compatibility}$. That is, $u^h$ is continuous across the interface, $u^{h,(\alpha)}(x) = u^{h,(\beta)}(x), \forall\, x \in \Gamma^{(\alpha,\beta)}$.

---

See Figure 6. Detailed implementation steps of the $G^0$ compatibility on B-splines and NURBS are provided in the next section.

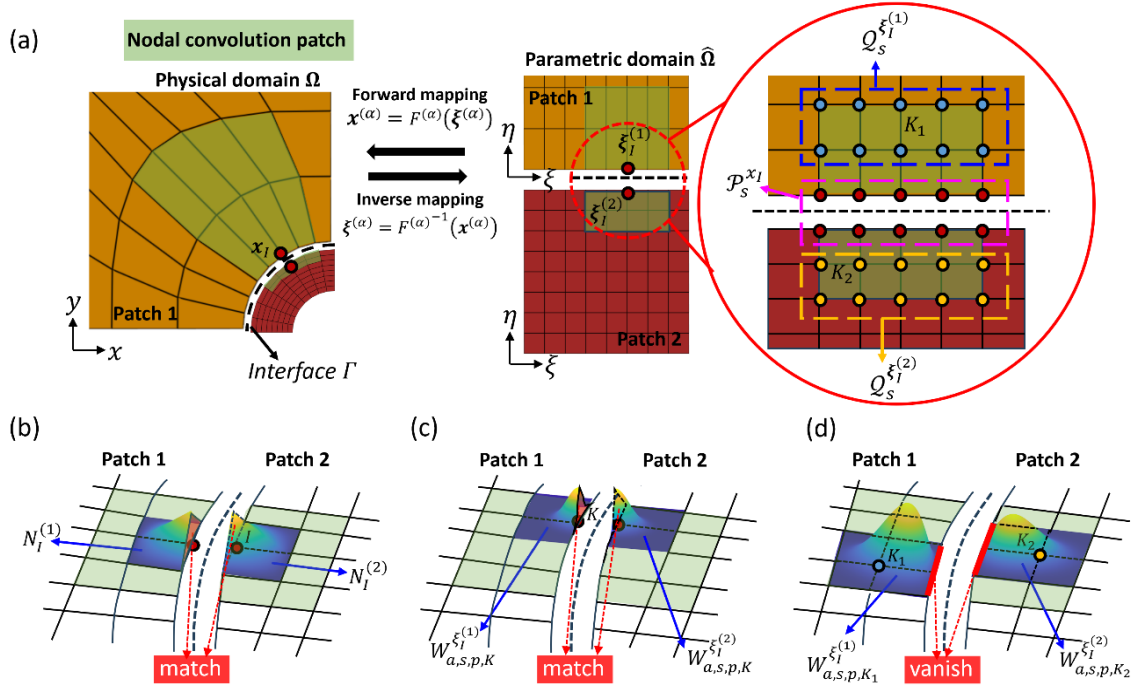

Figure 6 Illustration for compatibility conditions. (a) Physical space and parametric space. The nodal

convolution patches $A_s^{\xi_I^{(1)}}$, $A_s^{\xi_I^{(2)}}$ for shared node $I$ are colored by green. $A_s^{\xi_I^{(1)}}$, $A_s^{\xi_I^{(2)}}$ are divided into the index sets $Q_s^{\xi_I^{(1)}}$, $Q_s^{\xi_I^{(2)}}$ of internal nodes and the index set $\mathcal{P}_s^{x_I}$ of shared nodes along the interface. (b) Illustration for condition (C1): FEM shape functions $N_I^{(1)}, N_I^{(2)}$ for shared node $I$ match at the interface. (c) Illustration for condition (C2): convolution patch functions $W_{s,a,p,K}^{\xi_I^{(1)}}, W_{s,a,p,K}^{\xi_I^{(2)}}$ for node $K \in \mathcal{P}_s^{x_I}$ match at the interface. (d) Illustration for condition (C3): convolution patch functions $W_{s,a,p,K_1}^{\xi_I^{(1)}}, W_{s,a,p,K_2}^{\xi_I^{(2)}}$ for node $K_1 \in Q_s^{\xi_I^{(1)}}, K_2 \in Q_s^{\xi_I^{(2)}}$ vanish at the interface.

### 4.3 Implementation of $G^0$ Compatibility Condition via Product Rule

To satisfy the three conditions in Theorem 2, we construct FEM shape functions and convolution patch functions via product rule over a regular parametric mesh around the interface. Here, we consider B-spline and NURBS cases. To simplify the notations, two 2D neighboring patches, Patch 1 ($\Omega^{(1)}$) and patch 2 ($\Omega^{(2)}$), are considered for illustration, and $\Gamma^{(1,2)}$ is their interface. Without loss of generality, we assume that the interface of the two patches lies along the $\xi$-direction in the parametric domain $\widehat{\Omega}$ (where $\boldsymbol{\xi}^{(\alpha)} = (\xi^{(\alpha)}, \eta^{(\alpha)})$ denotes the vector of parametric coordinates in patch $\alpha$).

For each shared node $I$ at the interface $\Gamma^{(1,2)}$, the FEM shape functions $N_I^{(\alpha)}, \alpha = 1,2$, and convolution patch functions $W_{s,a,p,K}^{\xi_I^{(\alpha)}}$ are expressed as

$$N_I^{(\alpha)}(\boldsymbol{\xi}^{(\alpha)}) = N_i(\xi^{(\alpha)})N_j(\eta^{(\alpha)}), \tag{37}$$

$$W_{s,a,p,K}^{\xi_I^{(\alpha)}}(\boldsymbol{\xi}^{(\alpha)}) = \frac{\rho(\xi_K^{(\alpha)}, \eta_K^{(\alpha)})}{\rho(\xi^{(\alpha)}, \eta^{(\alpha)})} W_{s,a,p,m}^{\xi_i^{(\alpha)}}(\xi^{(\alpha)}) W_{s,a,p,n}^{\eta_j^{(\alpha)}}(\eta^{(\alpha)}), \alpha = 1,2. \tag{38}$$

Here, in the regular mesh in $\widehat{\Omega}$, the nodal index $I$ is represented by $(i,j)$ with $i,j$ representing the index in $\xi$-direction and $\eta$-direction, respectively. Correspondingly, the nodal index $K$ is represented by $(m,n)$. The $N_i$ and $N_j$ are 1D FEM shape functions in $\xi$-direction and $\eta$-direction, respectively. The $W_{s,a,p,m}^{\xi_i^{(\alpha)}}$ and $W_{s,a,p,n}^{\eta_j^{(\alpha)}}$ are 1D convolution patch functions, which satisfy Kronecker delta and reproducing properties. The $\rho(\xi^{(\alpha)}, \eta^{(\alpha)})$ is a factor coming from the weighting function $W^{(\alpha)}(\xi^{(\alpha)}, \eta^{(\alpha)})$ in NURBS, defined by:

$$\rho(\xi^{(\alpha)}, \eta^{(\alpha)}) = \begin{cases} 1, & \text{B-spline basis,} \\ \dfrac{W^{(\alpha)}(\xi^{(\alpha)}, \eta^{(\alpha)})}{W^{(\alpha)}|_{\Gamma^{(1,2)}}}, & \text{NURBS basis.} \end{cases} \tag{39}$$

Namely, the weighting function $W^{(\alpha)}(\xi^{(\alpha)}, \eta^{(\alpha)})$ in patch $\alpha$ is scaled by its trace at the interface. For B-spline basis functions used in the IGA mapping, $W(\xi^{(\alpha)}, \eta^{(\alpha)})$ can be regarded as a constant 1, and thus $\rho(\xi^{(\alpha)}, \eta^{(\alpha)}) = 1$. At node $K$, $\rho(\xi^{(\alpha)}, \eta^{(\alpha)})$ is normalized by $\rho(\xi_K^{(\alpha)}, \eta_K^{(\alpha)})$ to ensure Kronecker delta property of $W_{s,a,p,K}^{\xi_I^{(\alpha)}}$. Based on these, at the interface $\Gamma^{(1,2)}$, the FEM shape functions $N_I^{(\alpha)}(x)$ reduces to $N_i(\xi^{(\alpha)})$, and the convolution patch function $W_{s,a,p,K}^{\xi_I^{(\alpha)}}$ reduces to $W_{s,a,p,m}^{\xi_i^{(\alpha)}}(\xi^{(\alpha)})$. Therefore, if $N_i(\xi^{(\alpha)})$ and $W_{s,a,p,m}^{\xi_i^{(\alpha)}}(\xi^{(\alpha)}), \alpha = 1,2$ are the same at the interface, the three compatibility conditions in Theorem 2 are satisfied automatically. See Theorem 3.

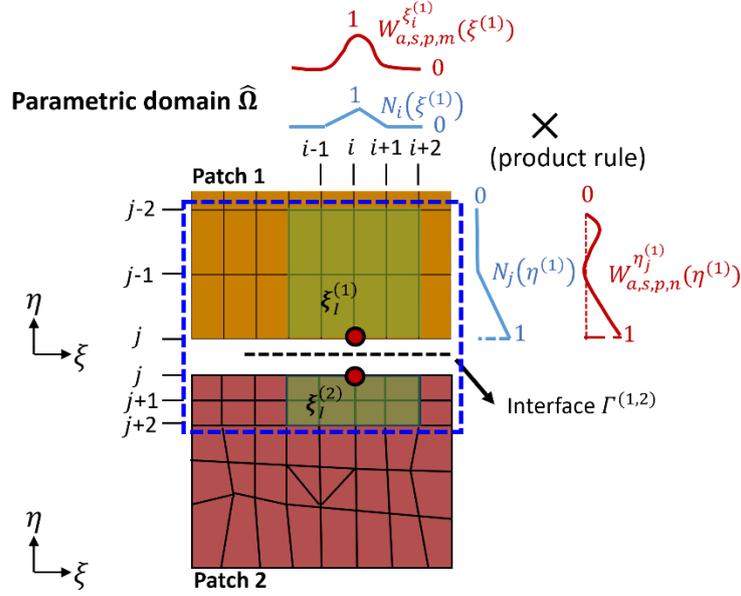

Figure 7 C-IGA shape functions are constructed by the product rule around the interface. In the region represented by a blue dashed box, regular meshes are established in the parametric domain $\hat{\Omega}$ where the FEM shape functions, and the convolution patch functions are constructed by the product rule defined in Eqs. ( 37 )-( 38 ). The nodal index $I$ is represented by $(i,j)$ with $i,j$ representing the index in $\xi$-direction and $\eta$-direction, respectively.

---

**Theorem 3**: **C-IGA basis constructed via product rule.** Let two neighboring patches $\Omega^{(1)}, \Omega^{(2)}$ satisfy the *nodal compatibility condition* with a patch interface $\Gamma^{(1,2)}$ that lies along the $\xi$-direction in the parametric domain $\hat{\Omega}$. Let C-IGA interpolations in two patches be constructed over regular parametric meshes in $\hat{\Omega}$ with the same parameters $s,a,p$.

Then for constructions of FEM shape functions $N_I^{(\alpha)}, \alpha = 1,2,$ and convolution patch functions $W_{s,a,p,K}^{\xi_I^{(\alpha)}}$ given in Eq. (37) and Eq. (38), if $N_i(\xi^{(\alpha)})$ and $W_{s,a,p,m}^{\xi_i^{(\alpha)}}(\xi^{(\alpha)}), \alpha = 1,2$ match at the interface, the three compatibility conditions in Theorem 2 are satisfied, and $\underline{C^0\ continuity}$ is attained at the interface $\Gamma^{(1,2)}$.

---

Next, we provide detailed constructions to ensure $N_i(\xi^{(\alpha)})$ and $W_{s,a,p,m}^{\xi_i^{(\alpha)}}(\xi^{(\alpha)}), \alpha = 1,2$ match at the interface. Node $I = (i,j)$ is at the interface.
- Two patches share the same $N_i(\xi^{(1)})$, i.e., $N_i(\xi^{(1)})$ in place of $N_i(\xi^{(2)})$ (or $N_i(\xi^{(2)})$ in place of $N_i(\xi^{(1)})$ ).
- Two patches share the same $W_{s,a,p,m}^{\xi_i^{(\alpha)}}(\xi^{(\alpha)})$ given by,

$$W_{s,a,p}^{\xi_i^{(1)}}(\xi^{(1)}) = W_{s,a,p,m}^{\xi_i^{(2)}}(\xi^{(2)}) = \Psi_a^T(x)A_\xi + P^T(\xi^{(1)},\xi^{(2)})K_\xi$$

with new basis

$$P(\xi^{(1)},\xi^{(2)}) = [P^{(1)}(\xi^{(1)}), P^{(2)}(\xi^{(2)})],$$

$$P^{(\alpha)}(\xi^{(\alpha)}) = \begin{cases} [1, (\xi^{(\alpha)}), \dots, (\xi^{(\alpha)})^p], & \text{B-spline basis,} \\ \dfrac{[1, (\xi^{(\alpha)}), \dots, (\xi^{(\alpha)})^p]}{W^{(\alpha)}|_{\Gamma^{(1,2)}}}, & \text{NURBS basis,} \end{cases} \quad \alpha = 1,2.$$

The $W_{s,a,p}^{\xi_i^{(\alpha)}}$ represents a vector of all convolution patch functions in the nodal convolution patch $A_s^{\xi_I^{(\alpha)}}$. The $\Psi_a(x)$ is a vector of radial basis functions defined in the physical domain. The $P(\xi^{(1)},\xi^{(2)})$ is a vector of basis functions to be reproduced, which are $p$-th order polynomial basis over a factor $W^{(\alpha)}|_{\Gamma^{(1,2)}}$. Coefficients $A_\xi$ and $K_\xi$ are determined to ensure Kronecker delta and reproducing properties. For B-spline basis functions, $P^{(\alpha)}(\xi^{(\alpha)})$ becomes a vector of polynomial basis.

In $\eta$-direction, $N_j(\eta^{(\alpha)})$ and $W_{s,a,p,m}^{\eta_j^{(\alpha)}}(\eta^{(\alpha)})$ are constructed in a similar way as 1D cases. Based on these constructions, the new convolution patch functions $W_{s,a,p,K}^{\xi_I^{(\alpha)}}$ preserve the Kronecker delta property and the ability to reproduce NURBS basis functions. The proof is provided in Appendix E.

Notably, the regular meshes and special constructions in Eq. (37, 38) are only required for elements that are adjacent to the interface within $s$ layers. For elements away from the interface, C-IGA shape functions are constructed similar to the single-patch C-IGA where irregular meshes are allowed.

## 5. Numerical Examples
### 5.1 1D Interpolations with Nodal Compatibility

We consider the following 1D interpolation problem to study the convergence of interface errors in a 2D multi-patch system. As illustrated in Figure 8, we adopt two different mappings $F^{(1)}, F^{(2)}$ over the same physical domain of interface $\Gamma^{(1,2)}: x \in [0,10]$ with the corresponding two different C-IGA interpolations $u^{(1)}(\xi^{(1)}), u^{(2)}(\xi^{(2)})$. The nodal compatibility suggests that $u^{(1)}(\xi^{(1)})$ and $u^{(2)}(\xi^{(2)})$ share the same physical nodes $[x_1, x_2, x_3, \dots]$ shown as red points in Figure 8. The corresponding nodal variables are denoted by $[u_1, u_2, u_3, \dots]$. Since different mappings for each patch lead to different constructions of shape functions, the discrepancy between the two C-IGA interpolations along the interface $\|u^{(1)}(x) - u^{(2)}(x)\|_{L^2(\Gamma^{(1,2)})}$ will not vanish.

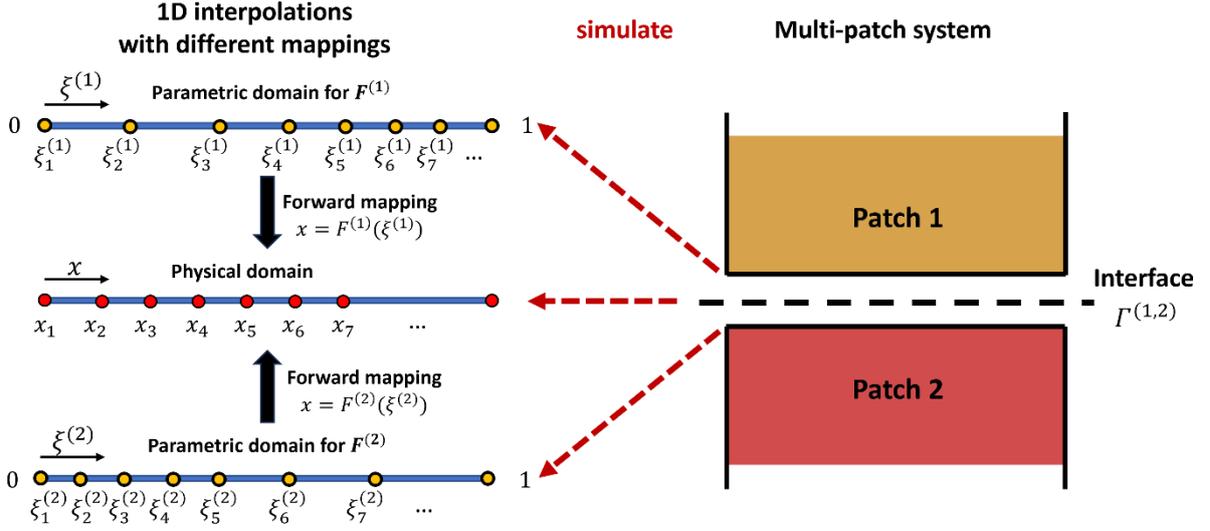

Figure 8 Illustration for 1D interpolation problem to simulate the interfaces of two patches. Two different mappings $F^{(1)}$ and $F^{(2)}$ generates two different interpolations $u^{(1)}$ and $u^{(2)}$, corresponding to interpolations in two neighboring patches. $u^{(1)}$ and $u^{(2)}$ share the same physical domain. The physical nodes are denoted by $[x_1, x_2, x_3, \ldots]$, associated with parametric nodes $[\xi_1^{(1)}, \xi_2^{(1)}, \xi_3^{(1)}, \ldots]$ for the first mapping $F^{(1)}$, and $[\xi_1^{(1)}, \xi_2^{(2)}, \xi_3^{(2)}, \ldots]$ for the second mapping $F^{(2)}$.

The two mappings are constructed with quadratic B-spline basis functions with control points $x = 0, 3, 10$, i.e., $x^{(1)} = F^{(1)}(\xi^{(1)}) = 0 \cdot (1 - \xi^{(1)})^2 + 3 \cdot 2\xi^{(1)}(1 - \xi^{(1)}) + 10 \cdot (\xi^{(1)})^2$ and $x = 0, 8, 10$, i.e., $x^{(2)} = F^{(2)}(\xi^{(2)}) = 0 \cdot (1 - \xi^{(2)})^2 + 8 \cdot 2\xi^{(2)}(1 - \xi^{(2)}) + 10 \cdot (\xi^{(2)})^2$, respectively. Two cases are considered: 1) when the nodal values are taken from a smooth 1D function, and 2) when the nodal values are taken from a discontinuous 1D function: $[u_1, u_2, u_3, u_4, \ldots] = [1, 2, 1, 2, \ldots]$. The former makes it as an interpolation of a smooth function, while the latter appears in the estimate of a solution of an elliptic PDE. The difference between the two interpolations is measured by the following relative $L^2$-norm deviation:

$$\frac{\left\|u^{(1)}(x) - u^{(2)}(x)\right\|_{L^2(\Gamma^{(1,2)})}}{\left\|u^{(1)}(x)\right\|_{L^2(\Gamma^{(1,2)})} + \left\|u^{(2)}(x)\right\|_{L^2(\Gamma^{(1,2)})}}. \quad (40)$$

First, let $u^{(1)}(x), u^{(2)}(x)$ be the interpolations of $u(x) = \sin x$, i.e., $u_i = \sin x_i$, $i = 1, 2, 3, \cdots$. The convergence of the relative $L^2$-norm deviation as the mesh size $h$ decreases is plotted in Figure 9(a). As expected from Lemma 1, the $L^2$ norm deviation reaches the $(p + 1)$-th order convergence rate. This can be understood from the following inequality:

$$\left\|u^{(1)}(x) - u^{(2)}\right\|_{L^2(\Gamma^{(1,2)})} \leq \left\|u^{(1)} - u\right\|_{L^2(\Gamma^{(1,2)})} + \left\|u^{(2)} - u\right\|_{L^2(\Gamma^{(1,2)})}. \quad (41)$$

C-IGA interpolation properties enable the $(p + 1)$-th order ($L^2$ norm) convergence rate for the two terms in the right-hand side. Thus, we have the $(p + 1)$-th order ($L^2$ norm) convergence rate for the interpolation of a smooth function in the multi-patch system.

Second, let the nodal values be $[u_1, u_2, u_3, u_4, ...] = [1,2,1,2, ...]$. Obviously, the frequency of oscillations increases as the mesh refines. As illustrated in Figure 9(b), the $L^2$ norm deviations exhibit an order 1 convergence rate rather than the $(p+1)$-th order convergence. In the error estimate of a PDE solution, the suboptimal convergence of the discontinuity across the interface obstructs the theoretical $p$-th order (for $H^1$ norm) convergence rate as indicated in Theorem 1.

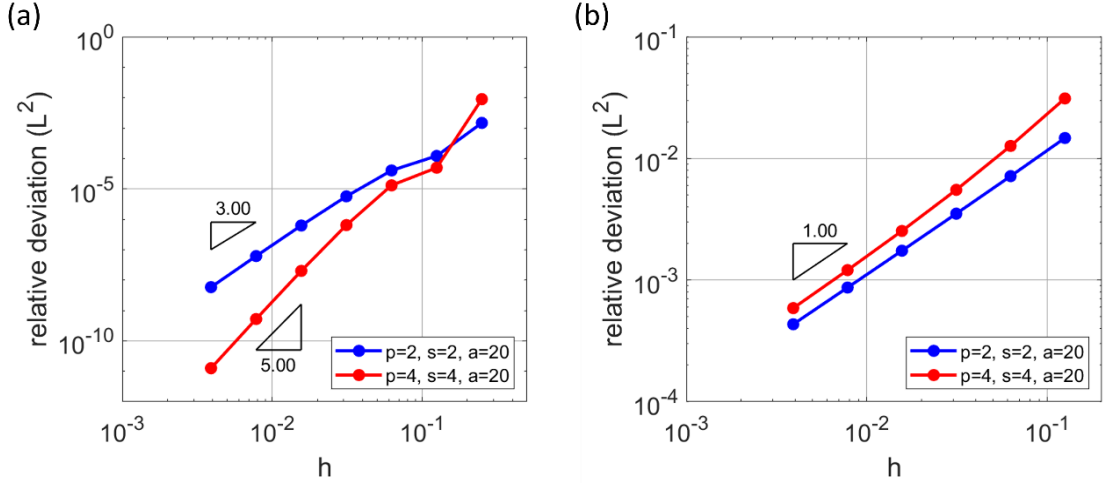

Figure 9 Convergence study of relative deviations between $u^{(1)}(x)$ and $u^{(2)}(x)$. (a) Case 1: $u^{(1)}(x)$ and $u^{(2)}(x)$ are interpolations of a smooth function, $u(x) = \sin x$. (b) Case 2: nodal values are assigned to be $[u_1, u_2, u_3, u_4, ...] = [1,2,1,2, ...]$. In Case 1, the $L^2$ norm relative deviation reaches the $(p+1)$-th order theoretical convergence rate. In Case 2, the convergence rate is only a first order.

### 5.2. 2D Interpolations with Nodal Compatibility

Next, we expand the problem definition to 2D and provide numerical proof on Lemma 1. Consider a plate with circular hole of radius $R = 0.5$, as shown in Figure 10(a). We interpolate a continuous function:

$$\sigma_{xx}(r, \theta) = 1 - \frac{R^2}{r^2}\left(\frac{3}{2}\cos 2\theta + \cos 4\theta\right) + \frac{3}{2}\frac{R^4}{r^4}\cos 4\theta, \quad (42)$$

which is known as an analytical stress field when a tensile force is applied to the plate [2]. Figure 10(a) shows the convergence of the $L^2$ norm error of the interpolation in the whole physical domain, $\Omega$. When $s = 2, p = 2$ are used for both patches, the convergence rate approaches to 3.0, which is the theoretical rate for $p = 2$. The same has been observed for $s = 3, p = 3$, which numerically proves Lemma 1.

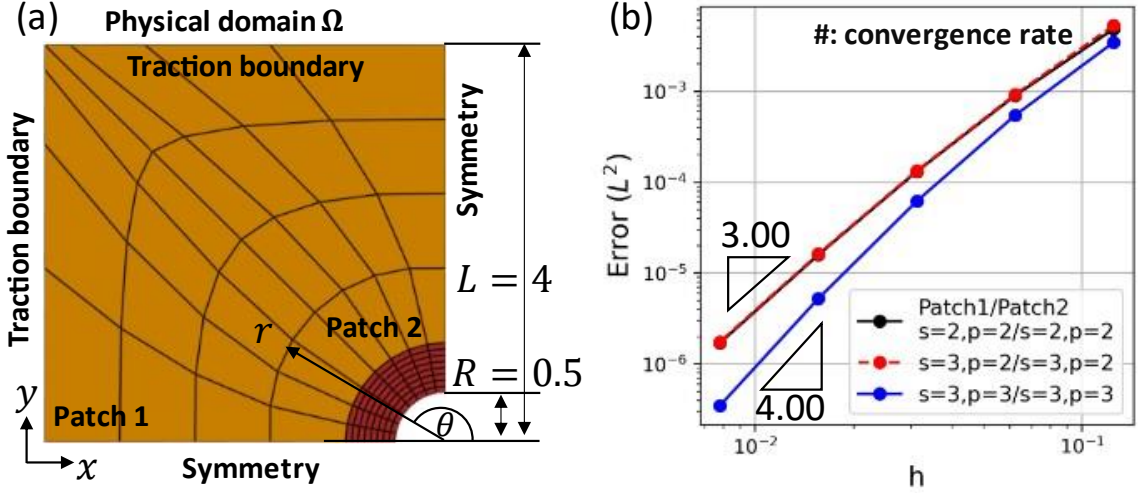

Figure 10 Convergence study of multi-patch C-IGA for interpolating a continuous function. (a) Geometrical description of a 2-patch problem; a plate with a circular hole. (b) Convergence rate with distinct polynomial order $(p)$ for each patch.

We have shown that multi-patch C-IGA with nodal compatibility can achieve theoretical convergence when it is used for interpolating a continuous function. However, as pointed out in Theorem 1, the theoretical convergence is not guaranteed for solving a PDE. We demonstrate this through an infinite plate with a circular hole problem [20, 34], schematically illustrated in Figure 10(a). The exact solution of this problem is given as:

$$\sigma_{xx}(r,\theta) = 1 - \frac{R^2}{r^2}\left(\frac{3}{2}\cos 2\theta + \cos 4\theta\right) + \frac{3}{2}\frac{R^4}{r^4}\cos 4\theta$$

$$\sigma_{yy}(r,\theta) = -\frac{R^2}{r^2}\left(\frac{1}{2}\cos 2\theta - \cos 4\theta\right) - \frac{3}{2}\frac{R^4}{r^4}\cos 4\theta, \quad (43)$$

$$\sigma_{xy}(r,\theta) = -\frac{R^2}{r^2}\left(\frac{1}{2}\sin 2\theta + \sin 4\theta\right) + \frac{3}{2}\frac{R^4}{r^4}\cos 4\theta.$$

where the right tractions are applied to the traction boundaries. As illustrated in Figure 10(a), the exact traction (Eq. (43)) is applied to the top and left end of the quarter domain under symmetric boundary condition and the hole radius is $R = 0.5$. The energy norm error is measured:

$$\|e_{energy}\| = \left(\int_\Omega \left(\frac{1}{2}\boldsymbol{C}_{mat}^{-1}(\boldsymbol{\sigma}^h - \boldsymbol{\sigma}_{exact}) \cdot (\boldsymbol{\sigma}^h - \boldsymbol{\sigma}_{exact})\right)d\Omega\right)^{1/2} \quad (44)$$

where $\boldsymbol{C}_{mat}^{-1}$ is the inverse of material stiffness matrix. Here, we assume a linear elasticity with the Young's modulus and Poisson's ratio being $10^3$ and 0.3, respectively. As shown in Figure 11(a), the solution convergence is sub-optimal for higher orders $p \geq 2$. This is attributed to the discontinuity at the interface; even though the nodal solutions match at the interface, the $G^0$ continuity is not achieved.

This is clearly shown in Figure 11(b), a plot for the $L^2$ norm deviation between the two patches at the interface (Eq. ( 40 )). The convergence rate remains constant as 1.0 for different polynomial orders $p$.

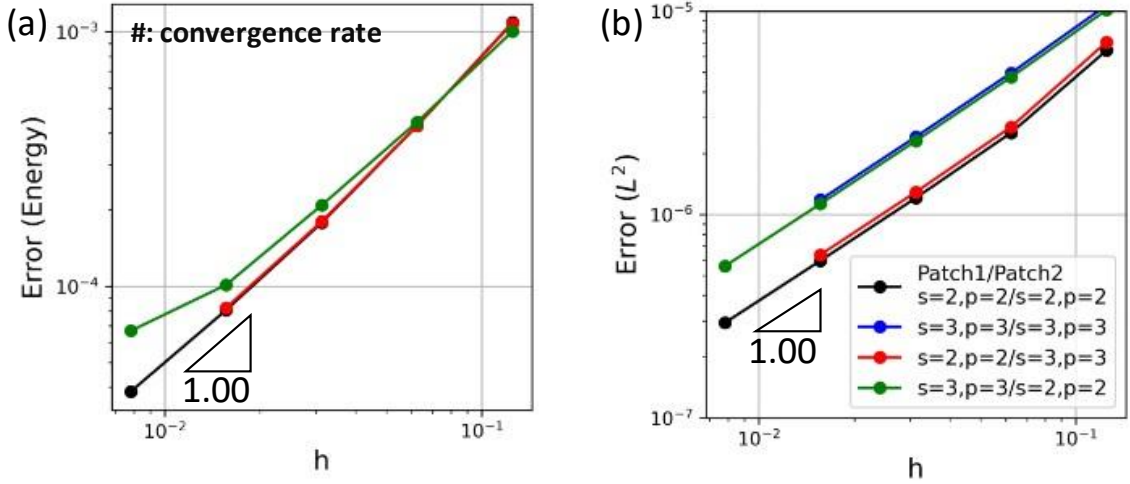

Figure 11 Convergence study of multi-patch C-IGA with nodal compatibility for solving a PDE (linear elasticity equation). (a) Energy norm error of the entire domain $\Omega$. (b) $L^2$ norm deviation at the interface.

### 5.3. Convergence Study of $G^0$ Compatibility

The $G^0$ compatibility condition is provided in Theorem 2, and detailed implementations using product rule are provided in Subsection 4.3 and Theorem 3. In this section, we provide numerical proofs on those compatibility conditions.

Consider a Poisson's equation within the same geometry $\Omega$ as in Figure 10(a), given by

$$\Delta u(x,y) + b(x,y) = 0, (x,y) \in \Omega. \qquad (45)$$

where $\Delta$ is the Laplace operator. We manufacture an exact solution as a Gaussian hump centered at $(x,y) = (-0.5, 1)$:

$$u^{Ext}(x,y) = e^{-\pi(x+0.5)^2 - \pi(x-1)^2}. \qquad (46)$$

The source term $b(x,y)$ is then derived from Eq. (42) using the exact solution. The boundary conditions are also determined from the exact solution, i.e.,

$$u(x,y)|_{\partial\Omega} = u^{Ext}(x,y)|_{\partial\Omega}. \qquad (47)$$

Figure 12 shows the geometry of the problem domain (a) and the exact solution (b). The interpolated solution field in these two patches are denoted by $u^{(1)}(x,y)$ and $u^{(2)}(x,y)$, respectively. The Gaussian hump locates near the interface denoted as the red dashed curve in Figure 12(a).

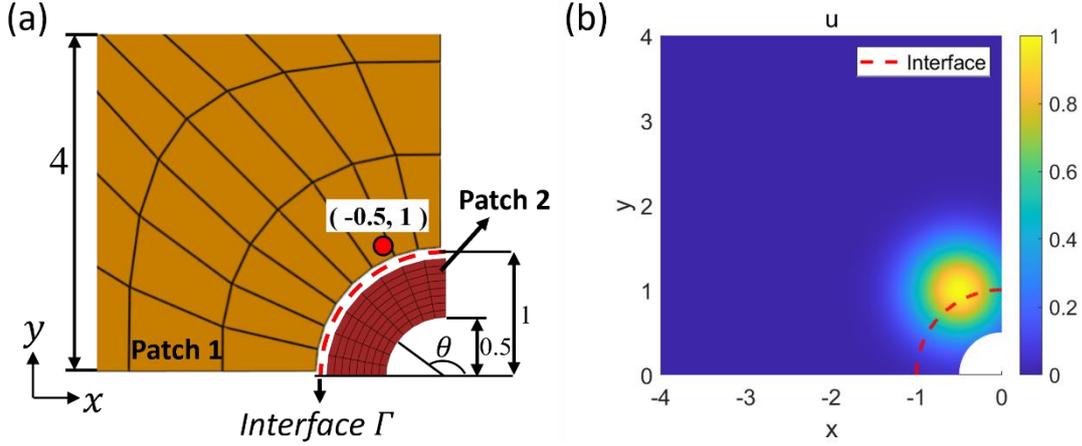

Figure 12 (a) Geometrical description of a 2-patch problem; a square plate with a circular hole. (b) Exact solution of the Poisson's equation. The red dashed curve represents the interface between the two patches.

It is worthwhile to note that in the coarsest IGA mesh (before any mesh refinement), as shown in Figure 13(a), two IGA elements exist in patch 1, while only one IGA element is in patch 2 [20]. When C-IGA shape functions are computed in patch 1, the convolution patch domains (the support domains of C-IGA shape functions) should not invade this IGA element boundary (same for patch 2) because the two patches have different parameterizations. If the convolution patch domains invade the IGA element boundary, the equivalence of geometric mappings among IGA and C-IGA will be violated.

For nodal compatibility, it is straightforward to construct C-IGA shape functions in each IGA element. However, if we want to achieve $G^0$ compatibility where the shape functions of the two patches coincide along the interface, it is required to divide patch 2 into two IGA elements similar to patch 1, as illustrated in Figure 13(b). Then the shape functions around the interface are constructed using the product rule. Detailed constructions are provided in Subsection 4.3, which ensures the three $G^0$ *compatibility conditions* in Theorem 2.

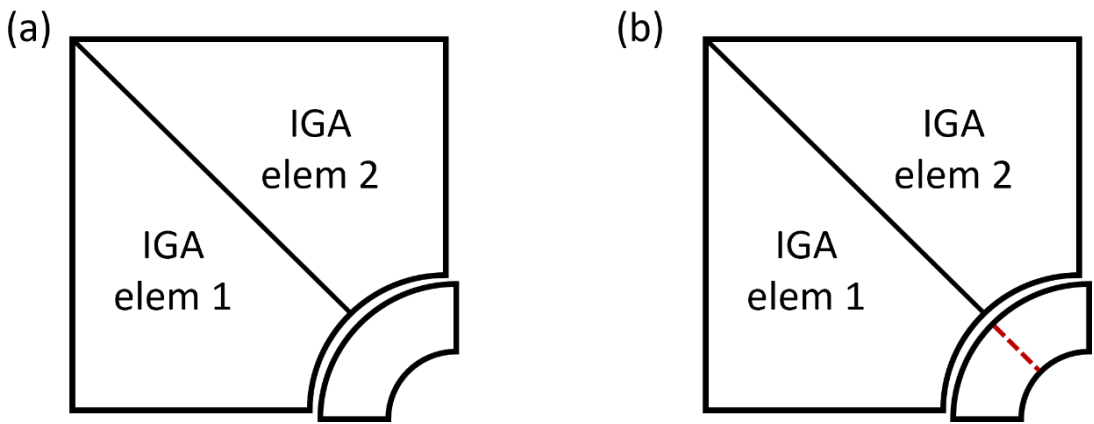

Figure 13 (a) In the coarsest mesh, two IGA elements are defined in patch 1 and one IGA element is needed in patch 2. (b) Patch 2 is divided into two parts for $G^0$ (b) Patch 2 is divided into two parts for $G^0$ compatibility.

In Figure 14(a), we plot $u^{(1)}(x,y)$ and $u^{(2)}(x,y)$ along the interface to evaluate the discontinuity; the two curves coincide with each other. As listed in Table 3, the relative deviation between $u^{(1)}(x,y)$ and $u^{(2)}(x,y)$ at the interface (defined in Eq. (37)) are less than 5e-12, which mostly originates from the machine precision. In contrast, the nodal compatibility leads to non-negligible relative deviations at the interface, varying from 1e-5 to 1e-7 depending on mesh refinement.

C-IGA shape functions with the $G^0$ compatibility at the interface is provided in Figure 15 *C-IGA shape functions with $G^0$ Figure 16*. The C-IGA shape functions encompass multiple C-IGA elements (not IGA elements) within a support domain of a normalized size $2(s+1) \times 2(s+1)$. They satisfy the Kronecker delta property, that is, for a given node, the corresponding shape function value becomes 1 at the node and 0 at the other nodes in the support domain of the shape function. As shown in Figure 15 *C-IGA shape functions with $G^0$ Figure 16*, C-IGA shape functions are allowed to have negative values within their support, which is different from IGA basis functions (i.e., splines).

We define the relative energy-norm error to measure the accuracy:

$$error = \frac{\|u^h - u^{Ext}\|_E}{\|u^{Ext}\|_E}, \|u\|_E = \int_\Omega |\nabla u|^2 \mathrm{d}x\mathrm{d}y. \quad (48)$$

The difference between the nodal compatibility and $G^0$ compatibility results in different convergence rates, as presented in Figure 14(b). The multi-patch C-IGA with $G^0$ compatibility can achieve theoretical convergence, i.e., $p$-th order energy-norm convergence rate. For $p=3, s=3$, C-IGA achieves a convergence rate of 3.61. However, the nodal compatibility results in a suboptimal convergence rate of 1.19 for $p=2, s=2$. As pointed out in Theorem 1, the solution accuracy is significantly influenced by the discontinuity across the interface.

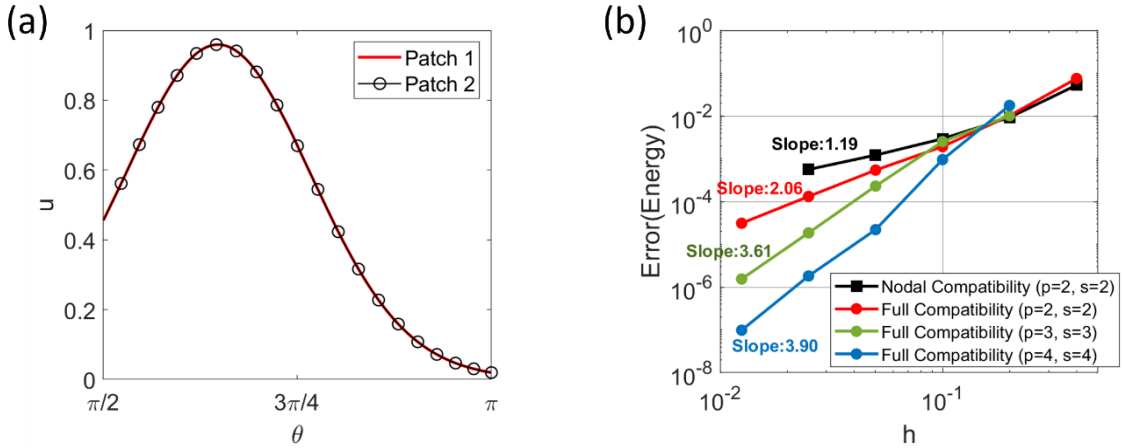

Figure 14 (a) Approximate solutions $u^{(1)}$ (patch 1) and $u^{(2)}$ (patch 2) along the interface in a polar coordinate system. (b) Convergence rate with distinct polynomial order $(p)$ for each patch. We set $s=p$ for every case.

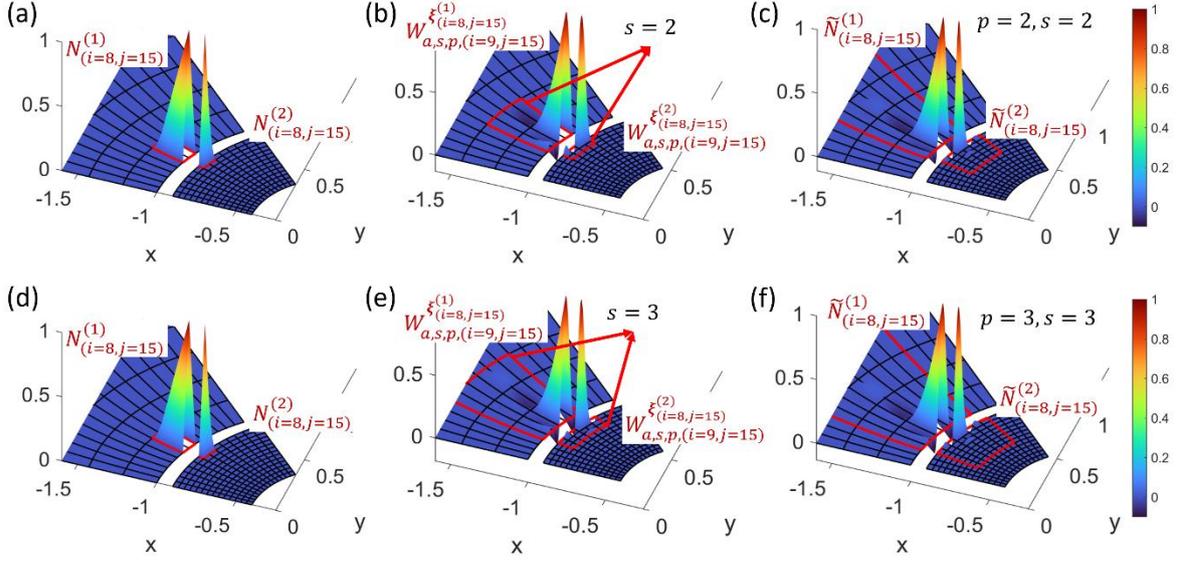

Figure 15 C-IGA shape functions with $G^0$   Figure 16 C-IGA shape functions with $G^0$ compatibility at the interface: (a) FEM shape functions $N^{(1)}_{(i=8,j=15)}$ and $N^{(2)}_{(i=8,j=15)}$ with $p=2, s=2$; (b) convolution patch functions $W^{\xi^{(1)}_{(i=8,j=15)}}_{a,s,p,(i=9,j=15)}$ and $W^{\xi^{(2)}_{(i=8,j=15)}}_{a,s,p,(i=9,j=15)}$ with $p=2, s=2$; (c) C-IGA shape functions $\tilde{N}^{(1)}_{(i=8,j=15)}$ and $\tilde{N}^{(2)}_{(i=8,j=15)}$ with $p=2, s=2$; (d) FEM shape functions $N^{(1)}_{(i=8,j=15)}$ and $N^{(2)}_{(i=8,j=15)}$ with $p=3, s=3$; (e) convolution patch functions $W^{\xi^{(1)}_{(i=8,j=15)}}_{a,s,p,(i=9,j=15)}$ and $W^{\xi^{(2)}_{(i=8,j=15)}}_{a,s,p,(i=9,j=15)}$ with $p=3, s=3$; (f) C-IGA shape functions $\tilde{N}^{(1)}_{(i=8,j=15)}$ and $\tilde{N}^{(2)}_{(i=8,j=15)}$ with $p=3, s=3$. There are $28 \times 28$ elements in the mesh.

Table 3 Relative deviations ($L^2$ norm defined by Eq. (37)) between $u^{(1)}(x)$ and $u^{(2)}(x)$ at the interface.

| Mesh | Nodal compatibility | | $G^0$ compatibility | |
|---|---|---|---|---|
| | $p=2, s=2$ | $p=3, s=3$ | $p=2, s=2$ | $p=3, s=3$ |
| $40 \times 40$ | 1.76E-05 | 2.42E-05 | 1.18E-12 | 3.82E-12 |
| $80 \times 80$ | 6.34E-06 | 9.91E-06 | 5.71E-13 | 4.03E-12 |
| $160 \times 160$ | 2.30E-06 | 3.32E-06 | 1.27E-12 | 4.64E-12 |
| $320 \times 320$ | 8.02E-07 | 1.16E-06 | 1.73E-12 | 1.98E-12 |
| $640 \times 640$ | - | - | 1.28E-12 | 4.63E-12 |

## 6. Conclusions and Future Works

Multi-CAD-patch (in short, multi-patch) C-IGA theory and its compatibility conditions are developed and validated with numerical examples. The multi-patch C-IGA formulation accompanies incompatibility of the solution field at patch interfaces, resulting in a suboptimal convergence of the solution accuracy when solving a partial differential equation (PDE). The two compatibility conditions are proposed and investigated: 1) nodal compatibility and 2) $G^0$ compatibility conditions. The former imposes the condition during the meshing stage, that is, the physical nodes of adjacent patches are set to be matched at the interface. This can be achieved by leveraging standard FEM meshing tools and inverse geometrical mapping of IGA. The nodal compatibility only guarantees $C^0$ continuity at those matching nodes. To achieve point-wise continuity across the interface, $G^0$ compatibility conditions are required. This is realized by constructing C-IGA shape functions with the product rule such that the shape functions at the interface are only affected by the nodes at the interface.

In terms of numerical implementation, the nodal compatibility is more practical because it only requires inverse mappings for each patch. The $G^0$ compatibility necessitates reparameterization of the near-interface region into a regular mesh. For most engineering applications where lots (even thousands) of CAD patches are involved for the geometrical representation, this reparameterization step will substantially impede the overall numerical analysis and complicate programming. In the future, the implementation of nodal compatibility to a realistic engineering problem with hundreds or thousands of CAD patches will be conducted.

Another critical and promising research direction is to leverage the multi-patch theory as a starting point of tensor decomposition (TD) – based reduced order modeling [16, 18, 35, 36]. In traditional TD or proper generalized decomposition (PGD) reduced order methods, a regular cartesian grid was indispensable, limiting their application to complex 3D geometries [16, 18, 35, 36]. The multi-patch IGA theories proposes a systematic reparameterization of a 3D complex CAD geometry into a regular parametric mesh [37, 38]. PGD has already been implemented to a single-patch domain in [39]. It is also recently reported that the C-HiDeNN-TD provides a powerful solver for space-time-parametric (S-T-P) PDEs where the spatial domain is represented with regular grid [32]. If is successfully merged with the multi-patch theory that provides a mapping from a regular mesh to a complex 3D geometry, it will be able to handle S-T-P PDEs with complex geometries, which has not been accomplished by other numerical methods.

**Acknowledgements**

The authors would like to thank Dr. Shaoqiang Tang, Dr. Dong Qian, and Dr. Gino Domel for their sincere and rigorous feedback.**Acknowledgements**

The authors would like to thank Dr. Shaoqiang Tang, Dr. Dong Qian, and Dr. Gino Domel for their sincere and rigorous feedback.

# Appendix A. Construction of convolution patch functions

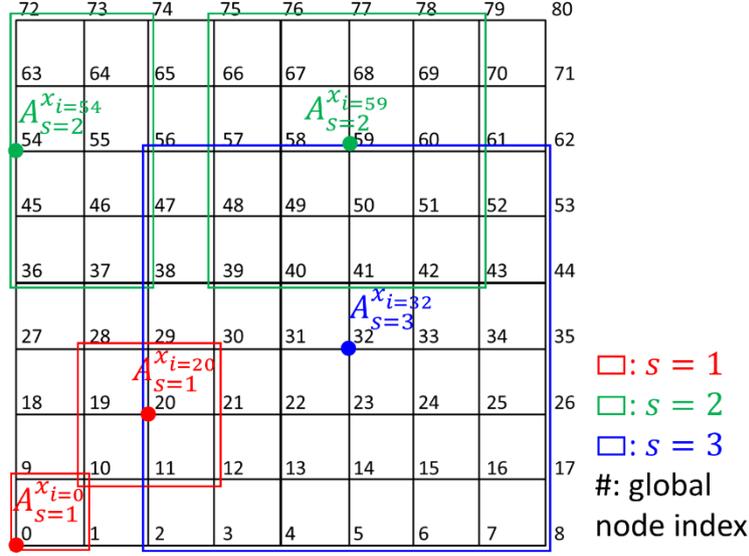

Figure 17 Illustration of the nodal convolution patch domain ($A_s^{x_i}$) on a two-dimensional (2D), 8 by 8 rectangular element mesh. The superscript $x_i$ denotes the center node of the nodal convolution patch, and the subscript is the patch size $s$. Each nodal convolution patch domain is colored by patch size $s$. The figure is borrowed from [17] with modification.

Convolution patch functions interpolate the nodal convolution patch domain $A_s^{x_i}$ where $s$ is the patch size (integer) and $x_i$ is the center node. Figure 17 illustrates how the nodal convolution patch domains are defined in a 2D mesh.

A field variable $u(x)$ inside a nodal convolution patch domain $A_s^{x_i}$ can be interpolated with radial basis $R_j(x)$ and polynomial basis $P_k(x)$ functions as

$$u(x) = \sum_{j=1}^{n} \Psi_j(x) A_j + \sum_{k=1}^{m} P_k(x) K_k \qquad (49)$$

$$= \Psi^T(x) A + P^T(x) K = [R^T(x) \quad P^T(x)] \begin{Bmatrix} A \\ K \end{Bmatrix}$$

where $n$ is the number of nodes in a nodal convolution patch domain $A_s^{x_i}$ and $m$ is the number of complete polynomial basis functions of a reproducing polynomial order $p$. For example, in 2D, the second order reproducing condition ($p = 2$) requires $m = 6$ polynomial basis functions: 1, $x$, $y$, $x^2, xy, y^2$. The $A$ and $K$ are the coefficients to be determined by imposing the Kronecker delta property.

A radial basis function $\Psi_j(x)$ is a function of radial distance between a point of interest $x$ and a nodal coordinate $x_j$. That is, $\Psi_j(x) = \Psi_j(r_j(x))$ where $r_j(x) = \|x - x_j\|_2$, the $L^2$ norm. Compactly supported radial basis functions listed in Table 3 are used, where the dimension of the local support is determined by the dilation parameter $a$. If $\frac{r_j}{a} > 1$, the radial basis functions return zero, making it compact-supported.

Table 3 Compactly supported radial basis functions $\Psi_j(r_j; a)$ used in [17] and this paper.

| Function | Formulation | Reference |
|---|---|---|
| Cubic spline | $\Psi_j(r_j; a) = \begin{cases} \dfrac{2}{3} - 4\dfrac{r_j^2}{a^2} + 4\dfrac{r_j^3}{a^3} & for\ 0 \leq \dfrac{r_j}{a} \leq \dfrac{1}{2} \\ \dfrac{4}{3} - 4\dfrac{r_j}{a} + 4\dfrac{r_j^2}{a^2} - \dfrac{4}{3}\dfrac{r_j^3}{a^3} & for\ \dfrac{1}{2} \leq \dfrac{r_j}{a} \leq 1 \\ 0 & otherwise \end{cases}$ | Liu (1995)[28] |
| Gaussian | $\Psi_j(r_j; a) = \begin{cases} \exp\left(-\dfrac{r_j^2}{a^2}\right) & for\ 0 \leq \dfrac{r_j}{a} \leq 1 \\ 0 & otherwise \end{cases}$ | |

The coefficients $A_j$ and $K_k$ are determined by imposing the Kronecker delta property, which can be written as

$$\boldsymbol{u}_n = \boldsymbol{\Psi}_n \boldsymbol{A} + \boldsymbol{P}_m \boldsymbol{K} \tag{50}$$

where the vector for the nodal field variable is

$$\boldsymbol{u}_n = \{u_1, u_2, \cdots u_n\}^T, \tag{51}$$

the moment matrix of the radial basis function is:

$$\boldsymbol{\Psi}_n = \begin{bmatrix} R_1(r_1) & R_2(r_1) & \cdots & R_n(r_1) \\ R_1(r_2) & R_2(r_2) & \cdots & R_n(r_2) \\ \vdots & \vdots & \ddots & \vdots \\ R_1(r_n) & R_2(r_n) & \cdots & R_n(r_n) \end{bmatrix}_{(n \times n)}, \tag{52}$$

the moment matrix of the polynomial basis function is:

$$\boldsymbol{P}_m = \begin{bmatrix} 1 & x_1 & y_1 & \cdots & p_m(\mathbf{x}_1) \\ 1 & x_2 & y_2 & \cdots & p_m(\mathbf{x}_2) \\ \vdots & \vdots & \vdots & \ddots & \vdots \\ 1 & x_n & y_n & \cdots & p_m(\mathbf{x}_n) \end{bmatrix}_{(n \times m)}, \tag{53}$$

and the coefficients $\boldsymbol{A}$ and $\boldsymbol{K}$ are:

$$\begin{aligned} \boldsymbol{A} &= \{A_1 \quad A_2 \quad \cdots \quad A_n\}^T \\ \boldsymbol{K} &= \{K_1 \quad K_2 \quad \cdots \quad K_m\}^T. \end{aligned} \tag{54}$$

As seen in Eq. ( 54 ), there are total $(n + m)$ unknowns but only $n$ equations exist in Eq. ( 50 ). Therefore, we need extra $m$ equations to fully determine the coefficients. The following additional $m$ constraints are usually used to resolve the ambiguity and to make the matrix equation symmetric [17, 40]:

$$\boldsymbol{P}_m^T \boldsymbol{A} = 0. \tag{55}$$

Finally, we get a system of $(n + m)$ linear equations:

$$\begin{Bmatrix} \boldsymbol{u}_n \\ \boldsymbol{0} \end{Bmatrix} = \begin{bmatrix} \boldsymbol{\Psi}_n & \boldsymbol{P}_m \\ \boldsymbol{P}_m^T & \boldsymbol{0} \end{bmatrix} \begin{Bmatrix} \boldsymbol{A} \\ \boldsymbol{K} \end{Bmatrix} = \boldsymbol{G} \begin{Bmatrix} \boldsymbol{A} \\ \boldsymbol{K} \end{Bmatrix}. \tag{56}$$

Note that the assembled moment matrix $\boldsymbol{G}$ is symmetric because the radial basis moment matrix $\boldsymbol{\Psi}_n$ is also symmetric and the additional $m$ constraints are set as the transpose of $\boldsymbol{P}_m$ (see Eq. ( 55 )). If

$G$ is a full-rank matrix, the coefficients $A$ and $K$ are uniquely determined as

$$\begin{Bmatrix} A \\ K \end{Bmatrix} = G^{-1} \begin{Bmatrix} u_n \\ 0 \end{Bmatrix}. \qquad (57)$$

By substituting Eq. (49) with Eq. (57), we get

$$u(x) = [\Psi^T(x) \quad P^T(x)] G^{-1} \begin{Bmatrix} u_n \\ 0 \end{Bmatrix}$$
$$= \widetilde{W}(x) \begin{Bmatrix} u_n \\ 0 \end{Bmatrix}. \qquad (58)$$

Here, $\widetilde{W}(x)$ is a $(n+m)$-component row vector and the last $m$ components are ignored because they are multiplied by zeros. If we denote the first $n$ components of $\widetilde{W}(x)$ by $W(x)$ and refer to *convolution patch functions*, finally the field variable is interpolated using radial basis functions with reproducing polynomial order $p$ and the Kronecker delta property:

$$u(x) = W(x) u_n \qquad (59)$$

Likewise, the derivative of the convolution patch function can be obtained by replacing the basis functions with their derivatives:

$$\frac{\partial u(x)}{\partial x} = \left[ \frac{\partial \Psi^T(x)}{\partial x} \quad \frac{\partial P^T(x)}{\partial x} \right] G^{-1} \begin{Bmatrix} u_n \\ 0 \end{Bmatrix}$$
$$= \frac{\partial \widetilde{W}(x)}{\partial x} \begin{Bmatrix} u_n \\ 0 \end{Bmatrix}. \qquad (60)$$

Since we can arbitrarily choose reproducing order $p$, C-HiDeNN can achieve arbitrary convergence rates without increasing the global DoFs, resulting in superior accuracy compared to linear FEM. In addition, due to the Kronecker delta property, we do not need special treatment for boundary conditions. For more discussions on the radial point interpolation method, readers may refer to [31, 40].

**Appendix B. Proof of Lemma 1.**
We prove Lemma 1 in this section.
**Proof.**
According to the definition of $\|\cdot\|_{H_b^k(\Omega)}$ in Eq. (17), we have

$$\|u - u^h\|_{H_b^1(\Omega)} = \left( \sum_{\alpha=1}^N \left( \|u^{(\alpha)} - u^{h,(\alpha)}\|_{H^1(\Omega^{(\alpha)})} \right)^2 \right)^{1/2} \leq \left( \sum_{\alpha=1}^N \left( C_{interp}^{(\alpha)} h^p \|u\|_{H^{p+1}(\Omega^{(\alpha)})} \right)^2 \right)^{1/2}$$
$$\leq \left( \sum_{\alpha=1}^N \left( C_{interp} h^p \|u\|_{H^{p+1}(\Omega^{(\alpha)})} \right)^2 \right)^{1/2} = C_{interp} h^p \left( \sum_{\alpha=1}^N \left( \|u\|_{H^{p+1}(\Omega^{(\alpha)})} \right)^2 \right)^{1/2}$$
$$= C_{interp} h^p \|u\|_{H^{p+1}(\Omega)}$$

with $C_{interp} = \max \left( C_{interp}^{(1)}, C_{interp}^{(2)}, \dots, C_{interp}^{(N)} \right)$.

## Appendix C. Proof of Theorem 1

We first define the following norms in the parametric domain $\widehat{\Omega}^{(\alpha)}$ for a single patch $\alpha$:

$$\|v\|^2_{H^1(\widehat{\Omega}^{(\alpha)})} = \int_{\widehat{\Omega}^{(\alpha)}} v^2 + |\nabla_\xi v|^2 dV(\xi^{(\alpha)}),$$

$$\|v\|^2_{L^2(\widehat{\Gamma}^{(\alpha)})} = \int_{\widehat{\Gamma}^{(\alpha)}} v^2 dV(\xi^{(\alpha)}), \quad v \in H^1(\Omega^{(\alpha)}),$$

where $\widehat{\Gamma}^{(\alpha)}$ is the boundary of $\widehat{\Omega}^{(\alpha)}$, $\nabla_\xi$ represents the derivatives with respect to parametric coordinates $\xi^{(\alpha)}$.

**Lemma 2:** Under the **Bijective Assumption,** the following norms in the physical domain $\Omega^{(\alpha)}$ and its corresponding parametric domain $\widehat{\Omega}^{(\alpha)} = [0,1]^d$ for a single patch $\alpha$ are equivalent, namely,

$$c^{(\alpha)}_{l1}\|v\|^2_{H^1(\widehat{\Omega}^{(\alpha)})} \leq \|v\|^2_{H^1(\Omega^{(\alpha)})} \leq c^{(\alpha)}_{l2}\|v\|^2_{H^1(\widehat{\Omega}^{(\alpha)})},$$

and

$$c^{(\alpha)}_{l3}\|v\|^2_{L^2(\widehat{\Gamma}^{(\alpha)})} \leq \|v\|^2_{L^2(\Gamma^{(\alpha)})} \leq c^{(\alpha)}_{l4}\|v\|^2_{L^2(\widehat{\Gamma}^{(\alpha)})}, \quad v \in H^1(\Omega^{(\alpha)}).$$

$c^{(\alpha)}_{l1}, c^{(\alpha)}_{l2}, c^{(\alpha)}_{l3}$ and $c^{(\alpha)}_{l4}$ are positive constants.

**Proof.**

Using the Jacobian transformation, we have:

$$\int_{\Omega^{(\alpha)}} v^2 dV(x) = \int_{\widehat{\Omega}^{(\alpha)}} v^2 \det(J^{(\alpha)}) dV(\xi^{(\alpha)}),$$

and

$$\int_{\Omega^{(\alpha)}} |\nabla v|^2 dV(x) = \int_{\widehat{\Omega}^{(\alpha)}} (\nabla_\xi v)^T \cdot (J^{(\alpha)})^{-T}(J^{(\alpha)})^{-1} \cdot \nabla_\xi v \det(J^{(\alpha)}) dV(\xi^{(\alpha)})$$

with $J^{(\alpha)}$ being the Jacobian matrix, i.e., $J^{(\alpha)} = \nabla(F^{(\alpha)})$. Thus, we have

$$\min|\det(J^{(\alpha)})| \int_{\widehat{\Omega}^{(\alpha)}} v^2 dV(\xi^{(\alpha)}) \leq \int_{\Omega^{(\alpha)}} v^2 dV(x) \leq \max|\det(J^{(\alpha)})| \int_{\widehat{\Omega}^{(\alpha)}} v^2 dV(\xi^{(\alpha)}),$$

and

$$\min|\det(J^{(\alpha)})| \cdot \lambda^{-1}_{max} \cdot \int_{\widehat{\Omega}^{(\alpha)}} |\nabla_\xi v|^2 dV(\xi^{(\alpha)}) \leq \int_{\Omega^{(\alpha)}} |\nabla v|^2 dV(x)$$

$$\leq \max|\det(J^{(\alpha)})| \cdot \lambda^{-1}_{min} \cdot \int_{\widehat{\Omega}^{(\alpha)}} |\nabla_\xi v|^2 dV(\xi^{(\alpha)}).$$

$\lambda_{min}$ and $\lambda_{max}$ are the minimum and maximum eigenvalues of $(J^{(\alpha)})^T J^{(\alpha)}$. According to the **Bijective Assumption**, we have

$$\min|\det(J^{(\alpha)})| > \delta > 0, (\lambda_{max})^{d-1}\lambda_{min} \geq \det\left((J^{(\alpha)})^T J^{(\alpha)}\right) > \delta^2 > 0.$$

Therefore, we have

$$(1 + \lambda^{-1}_{max})\min|\det(J^{(\alpha)})| \|v\|^2_{H^1(\widehat{\Omega}^{(\alpha)})} \leq \|v\|^2_{H^1(\Omega^{(\alpha)})} \leq (1 + \lambda^{-1}_{min})\max|\det(J^{(\alpha)})| \|v\|^2_{H^1(\widehat{\Omega}^{(\alpha)})}.$$

The second inequality is derived in the same manner.

**Lemma 3:** Under the **Bijective Assumption,** the following estimate holds:
$$\|v\|^2_{L^2(\Gamma^{(\alpha)})} \leq c_b^{(\alpha)} \|v\|^2_{H^1(\Omega^{(\alpha)})}, \qquad v \in H^1(\Omega^{(\alpha)}).$$

$c_b^{(\alpha)}$ is a positive constant.

**Proof.**

Note that the parametric domain $\hat{\Omega}^{(\alpha)}$ is a square or cubic domain $[0, 1]^d$. According to Gauss theorem, we have:
$$\frac{1}{2}\int_{\Gamma^{(\alpha)}} v^2 dV(\xi^{(\alpha)}) = \int_{\Gamma^{(\alpha)}} v^2 \left(\xi^{(\alpha)} - \xi_0^{(\alpha)}\right) \cdot \mathbf{n} dV(\xi^{(\alpha)}) = \int_{\hat{\Omega}^{(\alpha)}} \nabla \cdot \left(v^2 \left(\xi^{(\alpha)} - \xi_0^{(\alpha)}\right)\right) dV(\xi^{(\alpha)}),$$

where $\xi_0^{(\alpha)}$ is the center of $\hat{\Omega}^{(\alpha)}$, $\mathbf{n}$ the outer unit normal to $\Gamma^{(\alpha)}$, and thus $\left(\xi^{(\alpha)} - \xi_0^{(\alpha)}\right) \cdot \mathbf{n} = \frac{1}{2}$.

Moreover,
$$\int_{\hat{\Omega}^{(\alpha)}} \nabla \cdot \left(v^2 \left(\xi^{(\alpha)} - \xi_0^{(\alpha)}\right)\right) dV(\xi^{(\alpha)}) = \int_{\hat{\Omega}^{(\alpha)}} v^2 \nabla \cdot \left(\xi^{(\alpha)} - \xi_0^{(\alpha)}\right) + \left(\xi^{(\alpha)} - \xi_0^{(\alpha)}\right) \cdot \nabla v^2 dV(\xi^{(\alpha)})$$
$$= d\int_{\hat{\Omega}^{(\alpha)}} v^2 dV(\xi^{(\alpha)}) + 2\int_{\hat{\Omega}^{(\alpha)}} v\left(\xi^{(\alpha)} - \xi_0^{(\alpha)}\right) \cdot \nabla v dV(\xi^{(\alpha)})$$
$$\leq d\int_{\hat{\Omega}^{(\alpha)}} v^2 dV(\xi^{(\alpha)}) + 2\int_{\hat{\Omega}^{(\alpha)}} |v|\left|\xi^{(\alpha)} - \xi_0^{(\alpha)}\right||\nabla v| dV(\xi^{(\alpha)})$$
$$\leq d\int_{\hat{\Omega}^{(\alpha)}} v^2 dV(\xi^{(\alpha)}) + \frac{\sqrt{d}}{2}\int_{\hat{\Omega}^{(\alpha)}} \left(v^2 + |\nabla_\xi v|^2\right) dV(\xi^{(\alpha)})$$
$$\leq \left(d + \frac{\sqrt{d}}{2}\right)\int_{\hat{\Omega}^{(\alpha)}} \left(v^2 + |\nabla_\xi v|^2\right) dV(\xi^{(\alpha)}).$$

Thus, we have
$$\int_{\Gamma^{(\alpha)}} v^2 dV(\xi^{(\alpha)}) \leq \left(2d + \sqrt{d}\right)\int_{\hat{\Omega}^{(\alpha)}} \left(v^2 + |\nabla_\xi v|^2\right) dV(\xi^{(\alpha)}).$$

According to **Lemma 2, Lemma 4** is proven accordingly, namely,
$$\|v\|^2_{L^2(\Gamma^{(\alpha)})} \leq c_{l4}^{(\alpha)} \|v\|^2_{L^2(\hat{\Gamma}^{(\alpha)})} \leq c_{l4}^{(\alpha)}\left(2d + \sqrt{d}\right)\|v\|^2_{H^1(\hat{\Omega}^{(\alpha)})} \leq c_{l4}^{(\alpha)}\left(c_{l1}^{(\alpha)}\right)^{-1}\left(2d + \sqrt{d}\right)\|v\|^2_{H^1(\Omega^{(\alpha)})}.$$

Next the proof of **Theorem 1** is given as below:

**Proof.**

We define a new bilinear form including the penalty term:
$$A(w,v) = a(w,v) + \rho \sum_{\Gamma^{(\alpha,\beta)}\in\mathcal{F}} \left([w]_{\Gamma^{(\alpha,\beta)}}, [v]_{\Gamma^{(\alpha,\beta)}}\right)_{\Gamma^{(\alpha,\beta)}}, w,v \in H_b^1(\Omega).$$

According to Cauchy-Schwarz inequality, we have

$$A(w,v) \leq (a(w,w) \cdot a(v,v))^{\frac{1}{2}} + \rho \sum_{\Gamma^{(\alpha,\beta)}\in\mathcal{F}} \left(\|[w]\|_{L^2(\Gamma^{(\alpha,\beta)})} \cdot \|[v]\|_{L^2(\Gamma^{(\alpha,\beta)})}\right)$$

$$\leq \|w\|_{H_b^1(\Omega)} \cdot \|v\|_{H_b^1(\Omega)} + \rho \left(\sum_{\Gamma^{(\alpha,\beta)}\in\mathcal{F}} \|[w]\|_{L^2(\Gamma^{(\alpha,\beta)})}\right) \cdot \left(\sum_{\Gamma^{(\alpha,\beta)}\in\mathcal{F}} \|[v]\|_{L^2(\Gamma^{(\alpha,\beta)})}\right)$$

$$\leq \max(1,\rho) \left(\|w\|_{H_b^1(\Omega)} + \sum_{\Gamma^{(\alpha,\beta)}\in\mathcal{F}} \|[w]\|_{L^2(\Gamma^{(\alpha,\beta)})}\right) \cdot \left(\|v\|_{H_b^1(\Omega)} + \sum_{\Gamma^{(\alpha,\beta)}\in\mathcal{F}} \|[v]\|_{L^2(\Gamma^{(\alpha,\beta)})}\right)$$

Eqs. (3) and (5) imply that

$$A(w^h, u^h - u) = -\sum_{\Gamma^{(\alpha,\beta)}\in\mathcal{F}} \left([w^h]_{\Gamma^{(\alpha,\beta)}}, \langle\frac{\partial u}{\partial \boldsymbol{n}^{(\alpha)}}\rangle_{\Gamma^{(\alpha,\beta)}}\right)_{\Gamma^{(\alpha,\beta)}}, \forall w^h \in \mathcal{V}^h.$$

We denote $v^h$ as the interpolation of $u$. Obviously, $v^h \in \mathcal{S}^h$ and $u^h - v^h \in \mathcal{V}^h$. Then we obtain
$A(u^h - v^h, u^h - v^h) = A(u^h - v^h, u^h - u + u - v^h)$
$= A(u^h - v^h, u^h - u) + A(u^h - v^h, u - v^h)$

$$= -\sum_{\Gamma^{(\alpha,\beta)}\in\mathcal{F}} \left([u^h - v^h]_{\Gamma^{(\alpha,\beta)}}, \langle\frac{\partial u}{\partial \boldsymbol{n}^{(\alpha)}}\rangle_{\Gamma^{(\alpha,\beta)}}\right)_{\Gamma^{(\alpha,\beta)}} + A(u^h - v^h, u - v^h)$$

$$\leq \sum_{\Gamma^{(\alpha,\beta)}\in\mathcal{F}} \|[u^h - v^h]_{\Gamma^{(\alpha,\beta)}}\|_{L^2(\Gamma^{(\alpha,\beta)})} \cdot \left\|\langle\frac{\partial u}{\partial \boldsymbol{n}^{(\alpha)}}\rangle_{\Gamma^{(\alpha,\beta)}}\right\|_{L^2(\Gamma^{(\alpha,\beta)})}$$

$$+ \max(1,\rho) \left(\|u^h - v^h\|_{H_b^1(\Omega)} + \sum_{\Gamma^{(\alpha,\beta)}\in\mathcal{F}} \|[u^h - v^h]\|_{L^2(\Gamma^{(\alpha,\beta)})}\right)$$

$$\cdot \left(\|u - v^h\|_{H_b^1(\Omega)} + \sum_{\Gamma^{(\alpha,\beta)}\in\mathcal{F}} \|[u - v^h]\|_{L^2(\Gamma^{(\alpha,\beta)})}\right).$$

The estimate for the discontinuity implies that:

$$\|[u^h - v^h]_{\Gamma^{(\alpha,\beta)}}\|_{L^2(\Gamma^{(\alpha,\beta)})}$$

$$\leq C_{gap} \left(h^{\Gamma^{(\alpha,\beta)}}\right)^q \left(\|u^{h,(\alpha)} - v^{h,(\alpha)}\|_{L^2(\Gamma^{(\alpha,\beta)})} + \|u^{h,(\beta)} - v^{h,(\beta)}\|_{L^2(\Gamma^{(\alpha,\beta)})}\right).$$

Moreover, according to Lemma 3, we have

$$\|[u^h - v^h]_{\Gamma^{(\alpha,\beta)}}\|_{L^2(\Gamma^{(\alpha,\beta)})}$$

$$\leq C_{gap} \left(h^{\Gamma^{(\alpha,\beta)}}\right)^q \left(c_b^{(\alpha)}\|u^h - v^h\|_{H^1(\Omega^{(\alpha)})} + c_b^{(\beta)}\|u^h - v^h\|_{H^1(\Omega^{(\beta)})}\right)$$

$$\leq C_{gap}\left(h^{\Gamma^{(\alpha,\beta)}}\right)^q \left(c_b^{(\alpha)} + c_b^{(\beta)}\right) \|u^h - v^h\|_{H_b^1(\Omega)}.$$

According to **Lemma 1** and **Lemma 3,** we have

$$\|u - v^h\|_{H_b^1(\Omega)} + \sum_{\Gamma^{(\alpha,\beta)} \in \mathcal{F}} \|[u - v^h]\|_{L^2(\Gamma^{(\alpha,\beta)})}$$

$$\leq \|u - v^h\|_{H_b^1(\Omega)} + \sum_{\Gamma^{(\alpha,\beta)} \in \mathcal{F}} \left(\|u - v^h\|_{L^2(\Gamma^{(\alpha)})} + \|u - v^h\|_{L^2(\Gamma^{(\alpha)})}\right)$$

$$\leq \|u - v^h\|_{H_b^1(\Omega)} + \sum_{\Gamma^{(\alpha,\beta)} \in \mathcal{F}} \left(c_b^{(\alpha)} \|u - v^h\|_{H^1(\Omega^{(\alpha)})} + c_b^{(\beta)} \|u - v^h\|_{H^1(\Omega^{(\beta)})}\right)$$

$$\leq \left(1 + \sum_{\Gamma^{(\alpha,\beta)} \in \mathcal{F}} \left(c_b^{(\alpha)} + c_b^{(\beta)}\right)\right) \|u - v^h\|_{H_b^1(\Omega)}$$

$$\leq \left(1 + \sum_{\Gamma^{(\alpha,\beta)} \in \mathcal{F}} \left(c_b^{(\alpha)} + c_b^{(\beta)}\right)\right) C_{interp} h^p \|u\|_{H^{p+1}(\Omega)}.$$

In the same manner, we have

$$\|u^h - v^h\|_{H_b^1(\Omega)} + \sum_{\Gamma^{(\alpha,\beta)} \in \mathcal{F}} \|[u^h - v^h]\|_{L^2(\Gamma^{(\alpha,\beta)})} \leq \left(1 + \sum_{\Gamma^{(\alpha,\beta)} \in \mathcal{F}} \left(c_b^{(\alpha)} + c_b^{(\beta)}\right)\right) \|u^h - v^h\|_{H_b^1(\Omega)}$$

Thus, we have
$$A(u^h - v^h, u^h - v^h)$$

$$\leq \|u^h - v^h\|_{H_b^1(\Omega)} \cdot \left(\sum_{\Gamma^{(\alpha,\beta)} \in \mathcal{F}} \left(c_b^{(\alpha)} + c_b^{(\beta)}\right) C_{gap} \left(h^{\Gamma^{(\alpha,\beta)}}\right)^q \left\|\langle\frac{\partial u}{\partial \boldsymbol{n}^{(\alpha)}}\rangle_{\Gamma^{(\alpha,\beta)}}\right\|_{L^2(\Gamma^{(\alpha,\beta)})} + c_2 h^p \|u\|_{H^{p+1}(\Omega)}\right)$$

with $c_2 = \max(1, \rho) \left(1 + \sum_{\Gamma^{(\alpha,\beta)} \in \mathcal{F}} \left(c_b^{(\alpha)} + c_b^{(\beta)}\right)\right)^2 C_{interp}$, and

$$\|u^h - v^h\|_{H_b^1(\Omega)}^2 = a(u^h - v^h, u^h - v^h) \leq A(u^h - v^h, u^h - v^h).$$

This leads to

$$\|u^h - v^h\|_{H_b^1(\Omega)}$$

$$\leq \sum_{\Gamma^{(i,j)} \in \mathcal{F}} c_1^{\Gamma^{(\alpha,\beta)}} C_{gap} \left(h^{\Gamma^{(\alpha,\beta)}}\right)^\beta \left\|\langle\frac{\partial u^{Ext}}{\partial \boldsymbol{n}^{(\alpha)}}\rangle_{\Gamma^{(\alpha,\beta)}}\right\|_{L^2(\Gamma^{(\alpha,\beta)})} + c_2 h^p \|u^{Ext}\|_{H^{p+1}(\Omega)}$$

with $c_1^{\Gamma^{(\alpha,\beta)}} = c_b^{(\alpha)} + c_b^{(\beta)}$. Thus, the following error estimate holds:

$$\|u^h - u\|_{H_b^1(\Omega)} \leq \|u^h - v^h\|_{H_b^1(\Omega)} + \|u - v^h\|_{H_b^1(\Omega)}$$

$$\leq \sum_{\Gamma^{(\alpha,\beta)} \in \mathcal{F}} c_1^{\Gamma^{(\alpha,\beta)}} C_{gap} \left(h^{\Gamma^{(\alpha,\beta)}}\right)^\beta \left\|\langle\frac{\partial u^{Ext}}{\partial \boldsymbol{n}^{(\alpha)}}\rangle_{\Gamma^{(\alpha,\beta)}}\right\|_{L^2(\Gamma^{(\alpha,\beta)})} + (C_{interp} + c_2) h^p \|u^{Ext}\|_{H^{p+1}(\Omega)}.$$

## Appendix D. Proof of Theorem 2

**Theorem 2** provides three compatibility conditions for $G^0$ _compatibility_ The proof of **Theorem 2** is given as below:

**Proof.**

According to the Kronecker delta property of $N_I^{(\alpha)}$ and $W_{s,a,p,K}^{\xi_I^{(\alpha)}}$, we have

$$u^{(\alpha)}|_{\Gamma^{(\alpha,\beta)}} = \sum_{I|x_I \in \Gamma} N_I^{(\alpha)}|_{\Gamma^{(\alpha,\beta)}} \sum_{K \in A_s^{\xi_I^{(\alpha)}}} W_{s,a,p,K}^{\xi_I^{(\alpha)}}|_{\Gamma^{(\alpha,\beta)}} u_K + \sum_{I|x_I \notin \Gamma} N_I^{(\alpha)}|_{\Gamma^{(\alpha,\beta)}} \sum_{K \in A_s^{\xi_I^{(\alpha)}}} W_{s,a,p,K}^{\xi_I^{(\alpha)}}|_{\Gamma^{(\alpha,\beta)}} u_K$$

$$= \sum_{I|x_I \in \Gamma^{(\alpha,\beta)}} N_I^{(\alpha)}|_{\Gamma^{(\alpha,\beta)}} \sum_{K \in A_s^{\xi_I^{(\alpha)}}} W_{s,a,p,K}^{\xi_I^{(\alpha)}}|_{\Gamma^{(\alpha,\beta)}} u_K + 0 \cdot \sum_{K \in A_s^{\xi_I^{(\alpha)}}} W_{s,a,p,K}^{\xi_I^{(\alpha)}}|_{\Gamma^{(\alpha,\beta)}} u_K$$

$$= \sum_{I|x_I \in \Gamma} N_I^{(\alpha)}|_{\Gamma^{(\alpha,\beta)}} \left( \sum_{K \in A_s^{\xi_I^{(\alpha)}}, x_K \in \Gamma^{(\alpha,\beta)}} W_{s,a,p,K}^{\xi_I^{(\alpha)}}|_{\Gamma^{(\alpha,\beta)}} u_K + \sum_{K \in A_s^{\xi_I^{(\alpha)}}, x_K \notin \Gamma^{(\alpha,\beta)}} W_{s,a,p,K}^{\xi_I^{(\alpha)}}|_{\Gamma^{(\alpha,\beta)}} u_K \right)$$

$$= \sum_{I|x_I \in \Gamma^{(\alpha,\beta)}} N_I^{(\alpha)}|_{\Gamma^{(\alpha,\beta)}} \left( \sum_{K \in A_s^{\xi_I^{(\alpha)}}, x_K \in \Gamma^{(\alpha,\beta)}} W_{s,a,p,K}^{\xi_I^{(\alpha)}}|_{\Gamma^{(\alpha,\beta)}} u_K + 0 \right)$$

At the interface $\Gamma^{(\alpha,\beta)}$, the two patches $\Omega^{(\alpha)}, \Omega^{(\beta)}$ share the same physical nodes $x_K$ and nodal variable $u_K$. According to condition 2) and 3), we have

$$u^{(\alpha)}|_{\Gamma^{(\alpha,\beta)}} = \sum_{I|x_I \in \Gamma^{(\alpha,\beta)}} N_I^{(\alpha)}|_{\Gamma^{(\alpha,\beta)}} \sum_{K \in A_s^{\xi_I^{(\alpha)}}, x_K \in \Gamma^{(\alpha,\beta)}} W_{s,a,p,K}^{\xi_I^{(\alpha)}}|_{\Gamma^{(\alpha,\beta)}} u_K$$

$$= \sum_{I|x_I \in \Gamma^{(\alpha,\beta)}} N_I^{(\beta)}|_{\Gamma^{(\alpha,\beta)}} \sum_{K \in A_s^{\xi_I^{(\beta)}}, x_K \in \Gamma^{(i,j)}} W_{s,a,p,K}^{\xi_I^{(\alpha)}}|_{\Gamma^{(\alpha,\beta)}} u_K = u^{(\beta)}|_{\Gamma^{(\alpha,\beta)}}.$$

# Appendix E. Proof of Kronecker delta and reproducing properties for new C-IGA shape functions based on product rule.

The new convolution patch functions $W_{s,a,p,K}^{\xi_I^{(\alpha)}}$ in (41) based on product rule have the following properties.

**Theorem 4**: **Kronecker delta and reproducing properties.** Let $W_{s,a,p,K=(m,n)}^{\xi_{I=(i,j)}^{(\alpha)}}(\boldsymbol{\xi}^{(\alpha)})$ be convolution patch functions defined in a regular mesh in patch $\alpha$, given by

$$W_{s,a,p,K}^{\xi_I^{(\alpha)}}(\boldsymbol{\xi}^{(\alpha)}) = \frac{\rho(\xi_K^{(\alpha)},\eta_K^{(\alpha)})}{\rho(\xi^{(\alpha)},\eta^{(\alpha)})} W_{s,a,p,m}^{\xi_i^{(\alpha)}}(\xi^{(\alpha)}) W_{s,a,p,n}^{\eta_j^{(\alpha)}}(\eta^{(\alpha)}), \alpha = 1,2 \text{ with}$$

$$\rho(\xi^{(\alpha)},\eta^{(\alpha)}) = \frac{W^{(\alpha)}(\xi^{(\alpha)},\eta^{(\alpha)})}{W^{(\alpha)}|_\Gamma}.$$

$W^{(\alpha)}$ is the weighting function of NURBS basis functions, and $\Gamma$ is a patch interface of patch $\alpha$ that is along $\xi$-direction. If $W_{s,a,p,m}^{\xi_i^{(\alpha)}}(\xi^{(\alpha)})$ is a 1D convolution patch function in $\xi$-direction satisfying Kronecker delta and reproducing conditions

$$\sum_{m \in A_s^{\xi_i^{(\alpha)}}} W_{s,a,p,m}^{\xi_i^{(\alpha)}}(\xi^{(\alpha)}) \cdot \frac{\left(\xi_m^{(\alpha)}\right)^q}{W^{(\alpha)}\left(\xi_m^{(\alpha)}\right)|_\Gamma} = \frac{\left(\xi^{(\alpha)}\right)^q}{W^{(\alpha)}|_\Gamma}, q = 0,1,\ldots,p,$$

and $W_{s,a,p,n}^{\eta_j^{(\alpha)}}(\eta^{(\alpha)})$ is a 1D convolution patch function in $\eta$-direction satisfying Kronecker delta and reproducing conditions

$$\sum_{n \in A_s^{\eta_j^{(\alpha)}}} W_{s,a,p,n}^{\eta_j^{(\alpha)}}(\eta^{(\alpha)}) \cdot \left(\eta_n^{(\alpha)}\right)^q = \left(\eta^{(\alpha)}\right)^q, \quad q = 0,1,\ldots,p,$$

the following properties hold:

1) $W_{s,a,p,K}^{\xi_I^{(\alpha)}}\left(\boldsymbol{\xi}_J^{(\alpha)}\right) = \delta_{KJ}, K, J \in A_s^{\xi_I^{(\alpha)}}$ (Kronecker delta).

2) $\sum_{K \in A_s^{\xi_I^{(\alpha)}}} W_{s,a,p,K}^{\xi_I^{(\alpha)}}(\boldsymbol{\xi}^{(\alpha)}) \cdot \frac{\left(\xi_m^{(\alpha)}\right)^{q_1} \left(\eta_n^{(\alpha)}\right)^{q_2}}{W^{(\alpha)}\left(\xi_m^{(\alpha)},\eta_n^{(\alpha)}\right)} = \frac{\left(\xi^{(\alpha)}\right)^{q_1} \left(\eta^{(\alpha)}\right)^{q_2}}{W^{(\alpha)}(\xi^{(\alpha)},\eta^{(\alpha)})}, q_1, q_2 = 0,1,\ldots,p$ (Reproducing conditions).

**Proof.**
We first verify Kronecker delta property:

$$W_{s,a,p,K=(m,n)}^{\xi_{I=(i,j)}^{(\alpha)}}\left(\boldsymbol{\xi}_K^{(\alpha)}\right) = \frac{\rho\left(\xi_K^{(\alpha)},\eta_K^{(\alpha)}\right)}{\rho\left(\xi_K^{(\alpha)},\eta_K^{(\alpha)}\right)} W_{s,a,p,m}^{\xi_i^{(\alpha)}}\left(\xi_m^{(\alpha)}\right) W_{s,a,p,n}^{\eta_j^{(\alpha)}}\left(\eta_n^{(\alpha)}\right) = 1.$$

$$W_{s,a,p,K=(m,n)}^{\xi_{I=(i,j)}^{(\alpha)}}\left(\xi_{J=(r,t)}^{(\alpha)}\right) = \frac{\rho\left(\xi_K^{(\alpha)},\eta_K^{(\alpha)}\right)}{\rho\left(\xi_J^{(\alpha)},\eta_J^{(\alpha)}\right)} W_{s,a,p,m}^{\xi_i^{(\alpha)}}\left(\xi_r^{(\alpha)}\right) W_{s,a,p,n}^{\eta_j^{(\alpha)}}\left(\eta_t^{(\alpha)}\right)$$

$$= \frac{\rho\left(\xi_K^{(\alpha)},\eta_K^{(\alpha)}\right)}{\rho\left(\xi_J^{(\alpha)},\eta_J^{(\alpha)}\right)} \delta_{mr}\delta_{nt} = 0 \ (K \neq J).$$

We next verify reproducing conditions:

$$\sum_{K \in A_s^{\xi_I^{(\alpha)}}} W_{s,a,p,K}^{\xi_I^{(\alpha)}}(\boldsymbol{\xi}^{(\alpha)}) \cdot \frac{\left(\xi_m^{(\alpha)}\right)^{q_1}\left(\eta_n^{(\alpha)}\right)^{q_2}}{W^{(\alpha)}\left(\xi_m^{(\alpha)},\eta_n^{(\alpha)}\right)}$$

$$= \sum_{K \in A_s^{\xi_I^{(\alpha)}}} \frac{\rho\left(\xi_K^{(\alpha)},\eta_K^{(\alpha)}\right)}{\rho(\xi^{(\alpha)},\eta^{(\alpha)})} W_{s,a,p,m}^{\xi_i^{(\alpha)}}(\xi^{(\alpha)}) W_{s,a,p,n}^{\eta_j^{(\alpha)}}(\eta^{(\alpha)}) \cdot \frac{\left(\xi_m^{(\alpha)}\right)^{q_1}\left(\eta_n^{(\alpha)}\right)^{q_2}}{W\left(\boldsymbol{\xi}_K^{(\alpha)}\right)}$$

$$= \frac{1}{W(\boldsymbol{\xi}^{(\alpha)})} \sum_{K \in A_s^{\xi_I^{(\alpha)}}} W_{s,a,p,m}^{\xi_i^{(\alpha)}}(\xi^{(\alpha)}) W_{s,a,p,n}^{\eta_j^{(\alpha)}}(\eta^{(\alpha)}) \cdot \frac{W^{(\alpha)}\big|_\Gamma}{W^{(\alpha)}\left(\xi_m^{(\alpha)}\right)\big|_\Gamma} \cdot \left(\xi_m^{(\alpha)}\right)^{q_1}\left(\eta_n^{(\alpha)}\right)^{q_2}$$

$$= \frac{W^{(\alpha)}\big|_\Gamma}{W(\boldsymbol{\xi}^{(\alpha)})} \left(\sum_{m \in A_s^{\xi_i^{(\alpha)}}} W_{s,a,p,m}^{\xi_i^{(\alpha)}}(\xi^{(\alpha)}) \cdot \frac{\left(\xi_m^{(\alpha)}\right)^{q_1}}{W^{(\alpha)}\left(\xi_m^{(\alpha)}\right)\big|_\Gamma}\right) \left(\sum_{n \in A_s^{\eta_j^{(\alpha)}}} W_{s,a,p,n}^{\eta_j^{(\alpha)}}(\eta^{(\alpha)}) \cdot \left(\eta_n^{(\alpha)}\right)^{q_2}\right)$$

$$= \frac{W^{(\alpha)}\big|_\Gamma}{W(\boldsymbol{\xi}^{(\alpha)})} \frac{\left(\xi^{(\alpha)}\right)^{q_1}}{W^{(\alpha)}\big|_\Gamma} \left(\eta^{(\alpha)}\right)^{q_2} = \frac{\left(\xi^{(\alpha)}\right)^{q_1}\left(\eta^{(\alpha)}\right)^{q_2}}{W^{(\alpha)}(\xi^{(\alpha)},\eta^{(\alpha)})}.$$